\documentclass[11pt]{article}
\usepackage[margin=1.08in]{geometry}
\usepackage{amssymb}
\usepackage{amsmath,amsthm}
\usepackage{mathtools}
\usepackage{graphicx}
\usepackage[utf8]{inputenc}
\usepackage[dvipsnames]{xcolor}
\usepackage{tikz}
\usetikzlibrary{arrows.meta,positioning,calc,decorations.markings}
\usepackage{tikz-cd}
\usepackage{cite}
\usepackage{hyperref}
\usepackage{enumitem}
\usepackage{thmtools}
\usepackage{orcidlink}

\usepackage[nameinlink,noabbrev]{cleveref}

\numberwithin{equation}{section}

\newtheorem{theorem}{Theorem}[section]
\newtheorem{proposition}[theorem]{Proposition}
\newtheorem{lemma}[theorem]{Lemma}
\newtheorem{corollary}[theorem]{Corollary}
\newtheorem{definition}[theorem]{Definition}
\newtheorem{remark}[theorem]{Remark}
\newcommand{\N}{{\mathbb{N}}}
\newcommand{\HH}{{\mathbb{H}}}
\newcommand{\FF}{{\mathbb{F}}}
\newcommand{\Z}{{\mathbb{Z}}}
\newcommand{\cO}{{\mathcal{O}}}
\DeclareMathOperator{\Isom}{Isom}
\DeclareMathOperator{\Link}{Link}
\newcommand{\monic}{\xhookrightarrow{}}
\newcommand{\restr}[1]{|_{#1}}
\newcommand{\abs}[1]{\left|#1\right|}
\newcommand{\ct}{Cannon--Thurston}
\newcommand{\cat}[1]{CAT($#1$)}

\tikzset{
 relation edge/.style={line width=.65pt,
  postaction={decorate},decoration={markings,
   mark=at position .52 with {\arrow{Stealth[length=4.4pt,width=3.9pt]}}}},
 relation wall/.style={MidnightBlue,line width=.9pt},
 relation label/.style={font=\footnotesize,inner sep=1.8pt}
}
\newcommand{\relationoctagon}[4]{
 \coordinate (r0) at (0,.85);
 \coordinate (r1) at ({.16*#1},1.65);
 \coordinate (r2) at ({.38*#1},2);
 \coordinate (r3) at ({.62*#1},2);
 \coordinate (r4) at ({.84*#1},1.65);
 \coordinate (r5) at (#1,.85);
 \fill[MidnightBlue!3] (0,0)--(#1,0)--(r5)--(r4)--(r3)--(r2)--(r1)--(r0)--cycle;
 \draw[relation edge] (0,0)--node[relation label,below] {$#2$} (#1,0);
 \foreach \letter [count=\j from 0] in {#3} {
  \pgfmathtruncatemacro{\next}{\j+1}
  \pgfmathsetmacro{\labelpos}{\j==0 ? .68 : (\j==4 ? .32 : .5)}
  \draw[relation edge] (r\j)--node[relation label,above,sloped,pos=\labelpos] {$\letter$} (r\next);
 }
 \node[font=\small,text=black!55] at ({#1/2},.8) {$#4$};
}
\newcommand{\relationwall}[4]{%
 \draw[relation edge,#4] (#1,0)--node[relation label,#3] {$#2$} (#1,.85);
}
\colorlet{pairalpha}{black}
\colorlet{pairbeta}{black}
\colorlet{pairgamma}{black}
\colorlet{pairdelta}{black}
\newcommand{\occurrenceoctagon}[3]{%
 \foreach \k in {0,...,7} {\coordinate (o\k) at ({247.5+45*\k}:#1);}
 \fill[MidnightBlue!3] (o0)--(o1)--(o2)--(o3)--(o4)--(o5)--(o6)--(o7)--cycle;
 \foreach \letter/\edgecolor [count=\k from 0] in {#2} {
  \pgfmathtruncatemacro{\next}{mod(\k+1,8)}
  \draw[relation edge,\edgecolor] (o\k)--(o\next);
  \node[font=\small,text=\edgecolor] at
   ({270+45*\k}:{#1*cos(22.5)+.38}) {$\letter$};
 }
 \node[font=\small] at (0,0) {$#3$};
}
\newcommand{\sidepair}[3]{%
 \draw[relation edge,#3] (0,.28)--node[above,font=\small] {$#1$} (2.3,.28);
 \draw[relation edge,#3] (2.3,-.28)--node[below,font=\small] {$#2$} (0,-.28);
 \draw[densely dotted,black!55] (0,.28)--(0,-.28) (2.3,.28)--(2.3,-.28);
}

\title{A surface-group action on real-hyperbolic space with no Cannon--Thurston map}

\author{Bruno Dular\thanks{The construction was found by GPT-6 Astra. See \Cref{sec:description_of_AI_use} for details on the use of AI.\vspace{10pt}
}\;
\orcidlink{0009-0006-7784-9944}\thanks{University of Luxembourg, FSTM, Department of Mathematics,
Maison du nombre, 6 Avenue de la Fonte,
L-4364 Esch-sur-Alzette, Luxembourg.
Email: \href{mailto:bruno.dular@uni.lu}{\texttt{bruno.dular@uni.lu}}}}

\date{November 2, 2026}

\begin{document}
\maketitle

\begin{abstract}
	We present discrete, faithful and purely loxodromic representations of a closed surface group of genus 431, and of a free group, into the group of isometries of a finite-dimensional real-hyperbolic space for which no Cannon--Thurston map exists. The construction combines three ingredients: the small-cancellation strategy employed by Baker--Riley to construct a subgroup of an ambient hyperbolic group with no Cannon--Thurston map, the realization of these groups by octagon complexes, and a theorem of Desgroseilliers--Haglund that supplies a convex cocompact faithful representation of the ambient hyperbolic group on a real-hyperbolic space.

    These examples admit \emph{ray} \ct\ maps, as defined by Halder, Mahan Mj and Sardar.
\end{abstract}

\setcounter{tocdepth}{1}
\tableofcontents

\newpage
\section{Introduction}\label{sec:introduction}

\subsection{Context and main results}\label{sec:main_result}

For $g\geq 2$, let $S_g$ be a closed orientable surface of genus $g$ and let $\Gamma_g\coloneqq\pi_1(S_g)$, whose Gromov boundary $\partial\Gamma_g$ is homeomorphic to the circle $S^1$. Given a proper Gromov hyperbolic space $X$, with basepoint $x$ and Gromov boundary $\partial X$, and a representation $\rho\colon\Gamma_g\to\Isom(X)$, does the orbit map $\Gamma_g\to X\colon\gamma\mapsto\rho(\gamma)x$ extend continuously to a map
\begin{equation*}
    \Gamma_g\cup\partial\Gamma_g\to X\cup\partial X,
\end{equation*}
that sends $\partial\Gamma_g$ to $\partial X$?
If it does, such an extension is called a \emph{\ct\ map} for the representation $\rho$, and it is automatically $\Gamma_g$-equivariant. Here, $\Gamma_g$ is endowed with a word metric relative to a finite generating set. This question can be asked with $\Gamma_g$ replaced by any hyperbolic group, and with $X$ a hyperbolic group with its word metric.

Mahan Mj~\cite{MR1622603} formulated the above \ct\ map existence problem in the general setting of hyperbolic groups and spaces, extending questions previously studied for Kleinian groups, most notably in the 1985 preprint of Cannon and Thurston, which later appeared in published form as~\cite{MR2326947}. We briefly recall the historical context below.

Convex cocompact representations always admit \ct\ maps, as their orbit maps are quasi-isometric embeddings and hence extend continuously to embeddings of the boundary; see~\cite[III.H.3.9]{MR1744486}.

The history of \ct\ maps starts with $X=\HH^3$; see \cite{MR3652816,MR3966793} and the references therein. Floyd showed that geometrically finite representations admit a \ct\ map~\cite{MR568933} (extended to all dimensions by Tukia~\cite{MR955838}). For geometrically infinite representations, the first positive result was obtained by Cannon and Thurston~\cite{MR2326947} for representations induced by the fiber of a closed hyperbolic $3$-manifold fibering over the circle. Minsky extended the result to geometrically infinite representations with \emph{bounded geometry}, i.e.\ with a uniform positive lower bound on translation lengths~\cite{MR1257060}. Mahan Mj proved the existence of \ct\ maps for all simply or doubly degenerate purely loxodromic representations~\cite{MR3126566} (and subsequently for finitely generated Kleinian groups~\cite{MR3652816}). The methods rely on the geometry of degenerate hyperbolic $3$-manifolds, understood by results of Minsky and Brock--Canary--Minsky~\cite{MR2630036,MR2925381}. Little is known when $X=\HH^n$, $n\geq 4$, outside of the geometrically finite case; see also~\cite{MR2402415}.

In the general case, where $X$ is not a real-hyperbolic space, Mahan Mj has established the existence of \ct\ maps for normal hyperbolic subgroups of hyperbolic groups and for trees of hyperbolic spaces~\cite{MR1604882,MR1622603}. However, Baker and Riley have shown that \ct\ maps do not always exist~\cite{MR3143716}. Matsuda and Oguni used their example to show that every non-elementary hyperbolic group, in particular every closed surface group, embeds in a hyperbolic group for which no \ct\ map exists~\cite{MR3176651}.

In \cite[Question 6.13]{MR3966793}, Mahan Mj asks the question whether, for all $n\in\N$, all discrete faithful purely loxodromic representations of $\Gamma_g$ on $\HH^n$ by isometries admit \ct\ maps. We answer this question in the negative.

The construction adapts the Baker--Riley example, realizing it as an octagon complex to which a theorem of Desgroseilliers and Haglund applies~\cite{MR3110583}. The latter theorem uses Agol's virtual specialness theorem~\cite{MR3104553} and generalizes a construction of Kapovich~\cite{MR2175160}.

\begin{restatable}{theorem}{maintheorem}\label{thm:main}
    There exist a positive integer $q_0$, a torsion-free hyperbolic group $K_0$, a faithful discrete convex cocompact representation
    \begin{equation*}
        \sigma_0\colon K_0\longrightarrow\Isom(\HH^{q_0}),
    \end{equation*}
    and an injection $\iota\colon\Gamma_{431}\hookrightarrow K_0$ such that $\rho=\sigma_0\circ\iota$ is purely loxodromic and not convex cocompact, and has no \ct\ map.  More precisely, there are elements $b,d\in\Gamma_{431}$ and a sequence $(h_n)_{n\in\N}$ in $\Gamma_{431}$ with
    \begin{equation}\label{eq:intro_witness}
        h_nd^+\longrightarrow b^+\quad\text{in }\partial\Gamma_{431},\qquad
        \rho(h_n)\rho(d)^+\longrightarrow\rho(b)^-\ne\rho(b)^+,
    \end{equation}
    where $\gamma^{\pm}$, resp.\ $\rho(\gamma)^{\pm}$, denotes the attracting/repelling endpoint of an element $\gamma\in\Gamma_{431}$ in $\partial\Gamma_{431}$, resp.\ of its action on $\HH^{q_0}$.

    Non-trivial elements of $\rho(\Gamma_{431})$ have translation lengths uniformly bounded away from zero.
\end{restatable}

\begin{remark}
    The same construction produces a discrete, faithful and purely loxodromic representation of the free group $\FF_{287}$ into $\Isom(\HH^{q_0})$ with no \ct\ map; see \Cref{cor:free_group_action_with_no_CT}.
\end{remark}

Adding handles gives counterexamples in every genus $g\geq 431$, with all the representations obtained from one common convex cocompact ambient group; see \Cref{sec:extension_to_large_genera}.

\begin{restatable}{theorem}{maintheoremhighergenera}\label{thm:all_genera}
    There exist a positive integer $q_*$, a torsion-free hyperbolic group $K_*$, and a faithful discrete convex cocompact representation $\sigma_*\colon K_*\to\Isom(\HH^{q_*})$ such that, for every $g\geq 431$, there is an embedded surface group $\Gamma'_g\leq K_*$, isomorphic to $\pi_1(S_g)$, for which the restriction $\rho_g=\sigma_*\restr{\Gamma'_g}$ is discrete, faithful and purely loxodromic.
    
    Each of these restrictions is not convex cocompact and admits no \ct\ map. They share a positive lower bound on non-zero translation lengths.
\end{restatable}

\subsection{Ray \ct\ maps}\label{sec:intro_ray_CT_maps}

Halder, Mahan Mj and Sardar~\cite{HMS2025rayCTmap} introduced the notion of \emph{ray \ct\ map} for a representation $\rho\colon\Gamma\to\Isom(X)$, which encodes the limiting endpoints (or \emph{landing points}) of images of geodesic rays without necessarily being a continuous extension of an orbit map. They develop criteria for the existence of a ray \ct\ map and show, in particular, that the double HNN variant of Baker--Riley's construction admits a ray \ct\ map~\cite[Theorem 6.11]{HMS2025rayCTmap}. Since the construction presented here is a variant of Baker and Riley's construction, we invoke their criterion to show that the representation from \Cref{thm:main} (as well as those in higher genera) also admits a ray \ct\ map; see~\Cref{sec:ray_CT_map}.

\begin{restatable}{proposition}{prop_ray_CT_maps}\label{prop:ray_CT_maps}
    The representations $\rho$ and $\rho_{431+k}$, for $k\geq 0$, from the above theorems, all admit ray \ct\ maps.
\end{restatable}

Part of the argument is the property that admitting a ray \ct\ map is stable under almost malnormal quasi-convex extensions of hyperbolic groups; see~\Cref{prop:ray_CT_map_stable_under_malnormal_quasi-convex_extensions}.

\subsection{Plan of the construction}

The starting point is the construction of Baker and Riley~\cite{MR3143716} of a hyperbolic group $G$ and a free subgroup $F\leq G$ of finite rank, such that the inclusion $F\leq G$ admits no \ct\ map $\partial F\to\partial G$, using small-cancellation groups~\cite{MR577064,MR642423}.

They remark that, upon modifying the construction, one can obtain a group inclusion $F\leq G$ with no \ct\ map where $G$ is a \cat{-1}\ group~\cite[Remark 9]{MR3143716}. This requires specific choices of relations in $G$ that can be realized as boundary words of polygonal complexes, following Wise~\cite{MR1423338} (see also Rips~\cite{MR642423}). Related polygonal constructions were used by Barnard, Brady and Dani~\cite{MR2308946} to obtain highly distorted free subgroups of \cat{-1}\ groups. Barnard and Brady~\cite{MR2252898} also used surface attachments to construct exponentially distorted closed surface subgroups of \cat{0}\ hyperbolic groups.

The construction presented in this paper adapts the aforementioned ones, realizing them as octagon complexes, in order to apply a theorem of Desgroseilliers and Haglund~\cite[Theorem 1.7]{MR3110583} which provides the existence of a convex cocompact representation of $G$ on a real-hyperbolic space.

The obstruction to the existence of a \ct\ map that is used here is presented in \Cref{sec:CT_map}. It is different from the one used by Baker and Riley. Given a purely loxodromic representation $\rho\colon \Gamma\to\Isom(X)$, there is a well-defined boundary map
\begin{equation}\label{eq:boundary_map_on_ends_of_loxodromics}
    \{\gamma^{\pm}\mid \gamma\in\Gamma-\{1\}\}\to\partial X\colon \gamma^{\pm}\mapsto\rho(\gamma)^{\pm},
\end{equation}
which must be continuous, with the subspace topology induced from $\partial\Gamma$, if a \ct\ map exists for $\rho$. Thus, it suffices to find an example where the above map is not continuous; see \Cref{lem:CT_map_obstruction}.

\Cref{sec:algebraic_construction} presents an iterated HNN extension $G$ and a free subgroup $F\leq G$, following Baker--Riley~\cite{MR3143716}. \Cref{sec:explicit_obstruction} shows that the inclusion $F\leq G$ yields the obstruction above. The argument starts as in~\cite[\S 6.4]{HMS2025rayCTmap}, where Halder, Mahan Mj and Sardar give an alternative argument for the non-existence of a \ct\ map for the Baker--Riley family of constructions, but uses the non-continuity of \eqref{eq:boundary_map_on_ends_of_loxodromics} to conclude.

In \Cref{sec:words_adapted_to_octagons}, this general construction specializes to specific alphabet sizes and choices of words, with controlled two-letter subwords~\cite{MR1423338}. The resulting group $G$ is fixed from there on, and is realized as an octagon complex $Z$ in \Cref{sec:rows_of_octagons}.

At this point, the free subgroup $F$ of $G$ already satisfies the obstruction above and, applying the theorem of Desgroseilliers--Haglund~\cite{MR3110583} to the universal cover of $Z$ (after replacing the faces by regular Euclidean octagons), this yields a representation of $F$ on a real-hyperbolic space with no \ct\ map. To obtain a closed surface group representation, a surface group $\Gamma_{431}$ is amalgamated to $G$ along $F$, and \Cref{sec:attachment_criterion,sec:gluing_the_surface} verify that this amalgamation can be realized as an octagon complex that still satisfies the conditions of the theorem of Desgroseilliers--Haglund. The proof of \Cref{thm:main} is then concluded in \Cref{sec:action_on_real-hyperbolic_space}.

\Cref{sec:extension_to_large_genera} describes how to extend the construction to all genera $g\geq 431$, by inserting handles in the octagon complexes.

Finally, \Cref{sec:ray_CT_map} discusses the existence of ray \ct\ maps for the above constructions, as defined by Halder, Mahan Mj and Sardar~\cite{HMS2025rayCTmap}.

\subsection{Description of AI use}\label{sec:description_of_AI_use}
The construction, including the specific combinatorial data (the words and gluing patterns), was found by the AI model GPT-6 Astra (OpenAI, October 2026). It was obtained as an answer to the following prompt: ``Try to prove or disprove continuity of the boundary map on endpoints of loxodromic elements for rank-one symmetric spaces, without assuming a Cannon--Thurston existence theorem.''  I asked this after unsuccessfully attempting to prove continuity of that map and discussing with the model the gaps in my tentative arguments, and whether the Baker--Riley example could be realized in a rank-one symmetric space.

It took 50 minutes for three agents working simultaneously to reach a counterexample and proof, with a surface of genus 10,811 and using Wise words~\cite{MR1423338}. Following further prompts asking the agents to reduce the genus using the same strategy, the counterexample and proof presented in this paper were obtained.

I wrote this exposition based on that construction and I checked every step of the argument. I assume responsibility for the content, but I do not claim credit for the discovery of the construction.

The TikZ code for the figures was generated by the AI model GPT-6.1 Sol, from my hand-drawn sketches.

\subsection*{Acknowledgements}
I am grateful to Fran{\c c}ois Labourie for introducing me to the question of the existence of \ct\ maps and for discussions that motivated this work. I also thank Mahan Mj for comments on an earlier draft and for suggesting that I consider ray \ct\ maps, which led to the addition of~\Cref{sec:ray_CT_map}, and Ulysse Remfort-Aurat for helpful comments. I was supported by the Luxembourg National Research Fund (FNR), under its Aides \`a la Formation-Recherche (AFR) program, grant 18890152, during the preparation of this paper.

\section{Cannon--Thurston map obstruction}\label{sec:CT_map}

First, we recall some background on hyperbolic groups; see~\cite[Part III]{MR1744486}. Unless stated otherwise, a hyperbolic group $H$ is endowed with a word metric relative to a choice of finite generating set. Any two such choices yield quasi-isometric word metrics, and we can omit this choice when considering \emph{coarse} notions.

\begin{definition}[\!\!{\cite[III.$\Gamma$.3.4]{MR1744486}}]\label{def:quasi-convex}
    A subspace $C$ of a geodesic metric space $X$ is \emph{quasi-convex} if there exists a constant $k>0$ such that, for all $x,y\in C$, every geodesic segment joining $x$ to $y$ is contained in the $k$-neighborhood of $C$.

    A subgroup $F$ of a hyperbolic group $H$ is \emph{quasi-convex} if $F$ is quasi-convex in the Cayley graph of $H$.
\end{definition}

\begin{lemma}[\!\!{\cite[III.$\Gamma$.3.6]{MR1744486}}]\label{lem:quasi-convex_quasi-isometric}
    Let $F$ be a finitely generated subgroup of a hyperbolic group $H$.
    \begin{enumerate}[label=(\roman*)]
        \item Quasi-convexity of $F$ in $H$ is independent of the choice of finite generating set for $H$.
        \item $F\leq H$ is quasi-convex if and only if $F$ is quasi-isometrically embedded in $H$.
    \end{enumerate}
\end{lemma}

Let us recall the following orbit-map version of the definition of \ct\ maps for hyperbolic metric spaces, introduced by Mahan Mj~\cite[p.~136]{MR1622603}.

\begin{definition}[\ct\ map~\cite{MR1622603}]\label{def:CT_map}
    Let $\Gamma$ be a hyperbolic group, $X$ a proper Gromov hyperbolic space with basepoint $x\in X$, and $\rho\colon\Gamma\to\Isom(X)$ a representation. A \ct\ map for $\rho$ is a continuous extension of the orbit map $\Gamma\to X\colon\gamma\mapsto\rho(\gamma)x$ to a map
    \begin{equation*}
        \Gamma\cup\partial\Gamma\to X\cup\partial X
    \end{equation*}
    that sends $\partial\Gamma$ into $\partial X$. Existence of a \ct\ map, and its restriction to $\partial\Gamma$ when it exists, do not depend on the choice of basepoint or on the choice of word metric.
\end{definition}

Note that orbit maps and \ct\ maps, when they exist, are $\Gamma$-equivariant.

The attracting and repelling endpoints of a loxodromic isometry $h$ of a Gromov hyperbolic space $X$ are denoted $h^+$ and $h^-$, respectively. If $\Gamma$ is a hyperbolic group, then an element $\gamma$ is loxodromic if the map $\Gamma\to\Gamma\colon\alpha\mapsto\gamma\cdot\alpha$ is a loxodromic isometry, for any (equivalently, some) choice of word metric. In that case, we can write $\gamma^{\pm}\in\partial\Gamma$.

If $\gamma\in\Gamma$ has infinite order and $\rho(\gamma)$ is loxodromic, then any \ct\ map must send $\gamma^{\pm}$ to $\rho(\gamma)^\pm$. From this observation, the obstruction stated below follows, and shall be used in the proof of \Cref{thm:main}. See also~\cite[Lemma 2.10]{MR3488025}.

\begin{lemma}[Obstruction]\label{lem:CT_map_obstruction}
    Let $X$ be a proper Gromov hyperbolic space. Let $F$ be a hyperbolic group. Let $b,d\in F$ be elements of infinite order and let $(h_n)_{n\in\N}$ be a sequence in $F$. Suppose that a representation $\rho\colon F\to\Isom(X)$ sends $b,d$ to loxodromic isometries and satisfies
    \begin{equation*}
        h_nd^+\longrightarrow b^+,\qquad\rho(h_n)\rho(d)^+\longrightarrow z,
    \end{equation*}
    where $z\in\partial X$ is distinct from $\rho(b)^+$. Then $\rho$ admits no \ct\ map.
\end{lemma}
\begin{proof}
    Suppose that $\rho$ admits a \ct\ map and let $f\colon\partial F\to\partial X$ be its restriction to the boundary. Since $\rho(b)$ and $\rho(d)$ are loxodromic, we must have $f(b^+)=\rho(b)^+$ and $f(d^+)=\rho(d)^+$ by the above observation. Continuity and $F$-equivariance of $f$ yield
    \begin{equation*}
        \rho(h_n)\rho(d)^+=\rho(h_n)f(d^+)=f(h_nd^+)\longrightarrow f(b^+)=\rho(b)^+\neq z,
    \end{equation*}
    contradicting the hypotheses.
\end{proof}

In \Cref{sec:explicit_obstruction}, we shall use the following basic results as well. See e.g. \cite[\S 8.1]{MR919829} and \cite[\S 2]{MR1313451} for the first one, \cite[8.1.G]{MR919829} for the second one.

\begin{lemma}\label{lem:disjoint_fixed_points}
    Let $X$ be a proper geodesic Gromov hyperbolic space and let $G\leq\Isom(X)$ be a discrete torsion-free subgroup. Two loxodromic elements in $G$ with no common non-zero powers, i.e.\ elements $g,h\in G$ with $\langle g\rangle\cap\langle h\rangle=\{1\}$, have disjoint fixed-point sets.
\end{lemma}

\begin{lemma}[Attracting/repelling dynamics]\label{lem:moving_dynamics}
    Let $X$ be a proper Gromov hyperbolic space and let $a\in\Isom(X)$ be a loxodromic isometry. For any sequence $(x_n)_{n\in\N}$ in $X\cup\partial X$ and any sequence of positive integers $(m_n)_{n\in\N}$ with $m_n\to +\infty$, if $x_n\to x\neq a^-$, then $a^{m_n}x_n\to a^+$.
\end{lemma}

\section{Algebraic construction}\label{sec:algebraic_construction}

The construction starts with the double HNN extension of Baker--Riley~\cite[Lemma 5]{MR3143716}, presented in \Cref{sec:double_HNN,sec:explicit_obstruction}. In~\Cref{sec:words_adapted_to_octagons}, the sizes of alphabets and the choices of words are fixed to be well-suited for the realization of these groups as octagon complexes in \Cref{sec:geometric_realization}. We denote by $\FF_{m}$ the free group on $m$ generators.

\subsection{Double HNN extension}\label{sec:double_HNN}
Let $M,N\geq 1$. Let $C^0$ be the free group on the letters $c_0,\dots,c_{M-1}$ and $D^0$ be the free group on $d_0,\dots,d_{N-1}$. A word in a free group with a fixed basis is \emph{positive} if it is made only of positive powers of the letters of the basis. Choose non-empty positive words $C_i$ in the $c$-letters and $D_{ij}$ in the $d$-letters, referred to as \emph{positive substitution words}. Assume that
\begin{align*}
    \phi&\colon C^0\to C^0\colon c_i\mapsto C_i\\
    \psi_i&\colon D^0\to D^0\colon d_j\mapsto D_{ij}\qquad (i=0,\dots,M-1)
\end{align*}
are injective. Consider the group
\begin{equation}
    G\coloneqq\Big\langle b,c_i,d_j\mid b^{-1}c_ib=C_i,\; c_i^{-1}d_jc_i=D_{ij}\Big\rangle
\end{equation}
where the generators and relations range over all $i$ and $j$.

Denote by $C\coloneqq\langle c_i\rangle$ and $D\coloneqq\langle d_j\rangle$ the subgroups of $G$ obtained as the images of $C^0$ and $D^0$, respectively. Also let $F\coloneqq \langle b,D\rangle$.

\begin{lemma}\label{lem:double_HNN}
    The subgroups $C=\langle c_i\rangle$ and $D=\langle d_j\rangle$ of $G$ are free on their generators, satisfy $C\cap D=\{1\}$, and
    \begin{equation*}
        F\coloneqq\langle b,D\rangle\cong\langle b\rangle\ast D\cong\FF_{N+1}.
    \end{equation*}
    The group $G$ is torsion-free.
\end{lemma}
\begin{proof}
    This follows from Britton's Lemma~\cite{MR168633} on HNN extensions; see also~\cite{MR32641,MR577064}. Indeed,
    \begin{equation*}
        J_0\coloneqq\Big\langle c_0,D^0\mid c_0^{-1}d_jc_0=D_{0j}\Big\rangle
    \end{equation*}
    is the HNN extension of $D^0$ by $\psi_0$, with stable letter $c_0$. Then let $J_i$ be the HNN extension of $J_{i-1}$ by $\psi_i$ with stable letter $c_i$, for each $i=1,\dots,M-1$. Finally, $G$ is the HNN extension of $J_{M-1}$ by $\phi$ with stable letter $b$. The statement then follows from Britton's Lemma, as in~\cite[Lemma 5]{MR3143716}. See also~\cite[IV, Theorem 2.1]{MR577064}.
\end{proof}

\subsection{Explicit obstruction}\label{sec:explicit_obstruction}
Let us explain how the algebraic construction from the preceding subsection prevents the existence of a \ct\ map. It uses specific choices of words, adapted from~\cite[Proof of Theorem 6.10]{HMS2025rayCTmap}, where Halder, Mahan Mj and Sardar give an alternative argument for the non-existence of a \ct\ map for the Baker--Riley construction, or variants thereof; see~\cite[\S 6.4]{HMS2025rayCTmap}. The obstruction exhibited here, however, is not the same as in the aforementioned paper.

Write $c\coloneqq c_0$ and $d\coloneqq d_0$. For $n\geq 1$, let
\begin{equation*}
    u_n\coloneqq b^{-n}cb^{n}=\phi^n(c),\qquad w_n\coloneqq u_n^{-1}du_n.
\end{equation*}
Since $u_n=\phi^n(c)$, $u_n$ can be written as a positive $c$-word. Thus $w_n$ is obtained from $d$ by iterated application of the $\psi_i$, hence it is also a positive $d$-word, non-empty since $d\neq 1$. Define
\begin{equation}\label{eq:explicit_obstruction_hn}
    h_n\coloneqq b^nw_{2n}=b^{-n}c^{-1}b^{2n}db^{-2n}cb^{2n}.
\end{equation}
In particular, $h_n=b^nw_{2n}\in \langle b,D\rangle=F$. In the free basis $b,d_0,\dots,d_{N-1}$ of $F$, the infinite reduced word $h_nd^{\infty}$ begins with $b^n$. Thus,
\begin{equation}\label{eq:before_explicit_obstruction}
    h_nd_F^+\longrightarrow b_F^+,
\end{equation}
where $d_F^+,b_F^+$ denote the attracting fixed points of $d$ and $b$ in $\partial F$, respectively.

\begin{proposition}\label{prop:explicit_obstruction}
    With the notation of this section, for every discrete faithful purely loxodromic representation $\sigma\colon G\to\Isom(\HH^q)$,
    \begin{equation}
        \sigma(h_n)\sigma(d)^+\longrightarrow\sigma(b)^-.
    \end{equation}
    Consequently, $\sigma\restr{F}$ has no \ct\ map.
\end{proposition}
\begin{proof}
    For clarity, we write $g_{\sigma}\coloneqq\sigma(g)$ for generators $g$ of $G$.
    
    The elements $b$ and $d$ generate a free subgroup of rank two, thus the same holds for the loxodromic isometries $b_{\sigma}$ and $d_{\sigma}$ of $\HH^q$ and they have no common non-zero powers. Also, $b_{\sigma}$ and $c_{\sigma}$ have no common non-zero powers (e.g.\ use that there is a well-defined homomorphism $G\to\Z$ with $b\mapsto 1$ and $c_i,d_j\mapsto 0$).

    By \Cref{lem:disjoint_fixed_points}, this implies that $\{b_{\sigma}^+,b_{\sigma}^-\}\cap\{d_{\sigma}^+,d_{\sigma}^-\}=\emptyset$ and that $\{b_{\sigma}^+,b_{\sigma}^-\}\cap\{c_{\sigma}^+,c_{\sigma}^-\}=\emptyset$. In particular,
    \begin{equation}\label{eq:distinct_endpoints_in_explicit_obstruction_prop}
        d_{\sigma}^+\neq b_{\sigma}^-,\qquad c_{\sigma}b_{\sigma}^+\neq b_{\sigma}^+,\qquad d_{\sigma}b_{\sigma}^-\neq b_{\sigma}^-,\qquad c_{\sigma}^{-1}b_{\sigma}^+\neq b_{\sigma}^+.
    \end{equation}
    By \eqref{eq:explicit_obstruction_hn}, $\sigma(h_n)=b_{\sigma}^n\sigma(w_{2n})=b_{\sigma}^{-n}c_{\sigma}^{-1}b_{\sigma}^{2n}d_{\sigma}b_{\sigma}^{-2n}c_{\sigma}b_{\sigma}^{2n}$. We determine the limit of $\sigma(h_n)d_{\sigma}^+$ by applying \Cref{lem:moving_dynamics} iteratively:
    \begin{align*}
        b_{\sigma}^{2n}d_{\sigma}^+ &\longrightarrow b_{\sigma}^+,\\
        b_{\sigma}^{-2n}c_{\sigma}b_{\sigma}^{2n}d_{\sigma}^+ &\longrightarrow b_{\sigma}^-,\\
        b_{\sigma}^{2n}d_{\sigma}b_{\sigma}^{-2n}c_{\sigma}b_{\sigma}^{2n}d_{\sigma}^+ &\longrightarrow b_{\sigma}^+,\\
        \sigma(h_n)d_{\sigma}^+=b_{\sigma}^{-n}c_{\sigma}^{-1}b_{\sigma}^{2n}d_{\sigma}b_{\sigma}^{-2n}c_{\sigma}b_{\sigma}^{2n}d_{\sigma}^+ &\longrightarrow b_{\sigma}^-.
    \end{align*}
    At each step, \Cref{lem:moving_dynamics} applies since the intervening isometry ($c_{\sigma}^{\pm 1}$ or $d_{\sigma}$) moves $b_{\sigma}^{+}$ away from $b_{\sigma}^{+}$, or $b_{\sigma}^{-}$ away from $b_{\sigma}^{-}$, by \eqref{eq:distinct_endpoints_in_explicit_obstruction_prop}.

    The last line and \eqref{eq:before_explicit_obstruction} yield, by \Cref{lem:CT_map_obstruction}, the conclusion that $\sigma\restr{F}$ admits no \ct\ map.
\end{proof}

\subsection{Words adapted to rows of octagons}\label{sec:words_adapted_to_octagons}

In this section, we shall fix the size of the alphabets and the choice of positive substitution words, with controlled two-letter subwords as in Wise's construction~\cite[\S 2]{MR1423338}, and in the work of Barnard, Brady and Dani~\cite{MR2308946}. These choices will be fixed for the rest of the construction, and are precisely suited to the geometric realization of $G$ as the fundamental group of an octagon complex, which will occupy the next section.

The numerical choices below are adapted to the rows of octagons constructed in \Cref{sec:rows_of_octagons}. The $C_i$-relations use three octagons and require $\abs{C_i}=17$, whereas the $D_{ij}$-relations use four octagons and require $\abs{D_{ij}}=23$. We insert two ``zero differences'' (between subscripts of consecutive letters, see $\Delta$ below) in each $C_i$-word and three in each $D_{ij}$-word, responsible for repeated-letter pairs across \emph{internal walls} (edges separating octagons in a single row). This leaves $14$ non-zero differences for the $C_i$-words, which we let be all the non-zero residues modulo $M=15$. For each $i$, the $D_{ij}$-words have $19$ non-zero differences: the $15$ corresponding blocks partition the $285$ non-zero residues modulo $N=286$. Varying the initial letter will then produce every ordered pair of distinct letters exactly once across each family, as verified below.

These choices also determine the genus of the surface attached in \Cref{sec:gluing_the_surface}. The free subgroup $F=\langle b,D\rangle$ has rank $T=N+1=287$, and the surface constructed in \Cref{sec:gluing_the_surface} has one vertex, $4T$ edges and $T$ octagons. Thus its genus is $g=(3T+1)/2=431$.

Fix
\begin{equation*}
    M=15,\qquad N=286.
\end{equation*}
Subscripts on $c$-letters should be read modulo 15, and those on $d$-letters modulo 286.

Consider the following tuples of non-negative integers: the 16-tuple
\begin{equation*}
    \Delta\coloneqq (1,2,3,4,0,5,6,7,8,0,9,10,11,12,13,14)
\end{equation*}
and, for each $i\in\{0,\dots,M-1\}$, the 22-tuple
\begin{align*}
    \begin{split}
        \Delta_i\coloneqq 19i\cdot (1,\dots,1) + (1,2,3,4,5,-19i,6,7,8,&9,-19i,10,11,\\ &12,13,-19i,14,15,16,17,18,19),
    \end{split}
\end{align*}
so that $\Delta_i$ has zeros in positions 6,11,16, and the values $19i+1,19i+2,\dots,19i+19$ in the remaining entries.

For each $i\in\{0,\dots,M-1\}$, let $C_i$ be the positive $c$-word of length 17 starting with $c_i$ and such that, for each $k\in\{1,\dots,16\}$, the difference between the subscripts of its $(k+1)$-st letter and its $k$-th letter equals the $k$-th entry of $\Delta$, modulo 15. Similarly, for each $i\in\{0,\dots,M-1\}$ and $j\in\{0,\dots,N-1\}$, let $D_{ij}$ be the positive $d$-word of length 23 starting with $d_j$ and whose differences of consecutive subscripts are given by $\Delta_i$, modulo $286$.

For example, the words $C_0$ and $D_{00}$ are
\begin{align}\label{eq:example_C0_D00}
    \begin{split}
        C_0&=c_0c_1c_3c_6c_{10}c_{10}c_{0}c_{6}c_{13}c_{6}c_{6}c_{0}c_{10}c_{6}c_{3}c_{1}c_{0},\\
        D_{00}&=d_{0}d_{1}d_{3}d_{6}d_{10}d_{15}d_{15}d_{21}d_{28}d_{36}d_{45}d_{45}d_{55}d_{66}d_{78}d_{91}d_{91}d_{105}d_{120}d_{136}d_{153}d_{171}d_{190},
    \end{split}
\end{align}
where the repeated letters $d_{15}$, $d_{45}$, $d_{91}$ come from the three zeros in $\Delta_0$.

\begin{lemma}\label{lem:words_injection}
    The following homomorphisms, defined on generators,
    \begin{align*}
        \phi&\colon C^0\to C^0\colon c_i\mapsto C_i,\\
        \psi_i&\colon D^0\to D^0\colon d_j\mapsto D_{ij}\qquad (i=0,\dots,M-1)
    \end{align*}
    are injective.

    Among the words $C_0,\dots,C_{M-1}$, every ordered pair of distinct $c$-letters occurs exactly once at a position of non-zero difference. The same holds for the words $D_{ij}$. Pairs consisting of the same letter repeated twice only appear at the zero-difference positions.
\end{lemma}
\begin{proof}
    To show injectivity, we shall apply Stallings' \emph{graph-immersion} criterion~\cite[Proposition 5.3]{MR695906}, as in~\cite[Proposition 2.2]{MR2308946}.

    Indeed, the map $\phi$ can be realized as the induced map between fundamental groups of a map of graphs from a bouquet of $M$ oriented loops, each loop subdivided into $17$ edges, into a bouquet of $M$ oriented loops, each being a single edge, and sending the $i$-th loop to the concatenation of $17$ loops given by the word $C_i$. Since $C_i$ starts and ends with $c_i$, for each $i$, all the words $C_i$ have pairwise distinct first letters and pairwise distinct last letters. Since they are also positive words, the map of graphs is an \emph{immersion}. Thus, by Stallings' criterion, $\phi$ is an injection.

    Now, fix $i$ and consider $\psi_i\colon D^0\to D^0$. By the same argument, it suffices to notice that all the words $D_{ij}$, $j\in\{0,\dots,N-1\}$, are positive and have pairwise distinct starting letters and pairwise distinct ending letters. Indeed, the starting letter of $D_{ij}$ is $d_j$, and its ending letter is the $d$-letter of subscript $j+\sum_{k=1}^{19}(19i+k)$ modulo $286$. Thus $\psi_i$ is injective by Stallings' criterion.

    For the $c$-words, fix a non-zero difference $r\in\{1,\dots,14\}$. It occurs at exactly one position in $\Delta$. If $s$ is the sum of the preceding entries, the corresponding pair in $C_i$ is $(c_{i+s},c_{i+s+r})$, with subscripts modulo $15$. As $i$ varies, $i+s$ runs through every residue exactly once. Since every ordered pair $(c_a,c_b)$ with $a\ne b$ has a unique difference $r=b-a$ modulo $15$, each such pair occurs exactly once among the $C_i$.

    For the $d$-words, the non-zero differences in the tuples $\Delta_i$ are the disjoint blocks $\{19i+1,\dots,19i+19\}$, which together partition $\{1,\dots,285\}$. Fix one such difference $r$, and let $t$ be the sum of the entries preceding its unique position in $\Delta_i$. The corresponding pair in $D_{ij}$ is $(d_{j+t},d_{j+t+r})$, with subscripts modulo $286$. As $j$ varies, $j+t$ runs through every residue exactly once. Thus every ordered pair of distinct $d$-letters occurs exactly once among the $D_{ij}$. The zero differences give precisely the repeated-letter pairs.

    The last assertion is immediate.
\end{proof}

\section{Geometric realization as an octagon complex}\label{sec:geometric_realization}

We use the polygonal realization strategy of Wise~\cite{MR1423338} (see also~\cite{MR2308946}) and the surface-attachment strategies of Barnard and Brady~\cite[\S 4]{MR2252898}, with explicit octagons and markings described below. More precisely, we will describe the pushouts
\[
    \begin{tikzcd}
        R \arrow[r, hook] \arrow[d, hook] & Z \arrow[d, hook] \\
        Y_0 \arrow[r, hook] & Q_0\coloneqq Z\cup_R Y_0,
    \end{tikzcd}
    \qquad\text{and}\qquad
    \begin{tikzcd}
        F \arrow[r, hook] \arrow[d, hook] & G \arrow[d, hook] \\
        \Gamma_{431} \arrow[r, hook] & K_0\coloneqq G\ast_F \Gamma_{431},
    \end{tikzcd}
\]
where $R$ is a bouquet of loops realizing $F$, $Z$ is a hyperbolic octagon complex realizing $G$, and $Y_0$ is a hyperbolic octagon surface complex realizing $\Gamma_{431}$.

\subsection{Octagon complexes}\label{sec:octagon_complexes}

Let us recall some notions of hyperbolic polygonal complexes. See~\cite[I.7]{MR1744486} for more details. Note that \cite{MR1744486} requires each polygonal cell of a polygonal complex to be embedded, while the ones considered here will have vertices and some sides identified. The results from \cite{MR1744486} still apply by passing to the \emph{barycentric} subdivision.

All octagons in this section are congruent regular right-angled hyperbolic octagons. Their sides are labelled by letters, with $a^{-1}$ denoting the side $a$ with opposite orientation. If an edge $a$ of an octagon complex is a loop, we distinguish its \emph{incoming} (or \emph{negative}) half-edge $[a]_-$ and \emph{outgoing} (or \emph{positive}) half-edge $[a]_+$. We have $[a^{-1}]_{\pm}=[a]_{\mp}$.

The \emph{(geometric) link} of a vertex $v$ in a hyperbolic polygonal complex is the metric graph whose vertices are half-edges incident to $v$, and with an edge of length $\theta$ joining two half-edges whenever they are consecutive sides of a polygon making angle $\theta$ at $v$. In particular, in the present construction, all link edges have length $\pi/2$.

We shall use the following criterion, as used by Wise~\cite{MR1423338}.

\begin{proposition}\label{prop:hyperbolicity_of_octagon_complex}
    Let $K$ be a finite connected complex whose faces are regular right-angled hyperbolic octagons, glued isometrically along edges. If its vertex links are simple, i.e.\ have neither loop edges nor parallel edges, and triangle-free, then $K$ is locally \cat{-1}\ and $\pi_1(K)$ is hyperbolic and torsion-free.
\end{proposition}
\begin{proof}
    Let $v$ be a vertex of $K$. Each edge of $\Link(v)$ has length $\pi/2$. Since $\Link(v)$ is simple and triangle-free, every embedded circuit in it has at least four edges and hence length at least $2\pi$. Thus $K$ has \cat{1}\ links (and, when passing to the barycentric subdivision, the same holds for the new vertices too), hence $K$ is locally \cat{-1}\ by the link condition~\cite[II.5.2, II.5.6]{MR1744486}. The universal cover $\widetilde{K}$ is simply connected and has the same local links, hence it is globally \cat{-1}\ by~\cite[II.5.4(4)]{MR1744486}. Also, finiteness of $K$ makes $\widetilde{K}$ proper and the action of $\pi_1(K)$ on $\widetilde{K}$ cocompact. Thus $\pi_1(K)$ is hyperbolic. Since it acts freely on $\widetilde{K}$, it is torsion-free~\cite[II.2.8(1)]{MR1744486}.
\end{proof}

The following result will be used several times, hence we recall it here.

\begin{proposition}[\!\!{\cite[II.4.14]{MR1744486}}]\label{prop:local_isometries_are_pi_1_injective}
    Let $X$ and $Y$ be complete connected metric spaces, such that $X$ is non-positively curved and $Y$ is locally a length space. If a map $f\colon Y\to X$ is a local isometric embedding, then $Y$ is also non-positively curved, the induced homomorphism $f_*\colon\pi_1(Y,y)\to\pi_1(X,f(y))$ is injective for any $y\in Y$ and any continuous lift of $f$ to the universal covers $\widetilde{f}\colon \widetilde{Y}\to\widetilde{X}$ is an isometric embedding.
\end{proposition}

\subsection{Rows of octagons}\label{sec:rows_of_octagons}

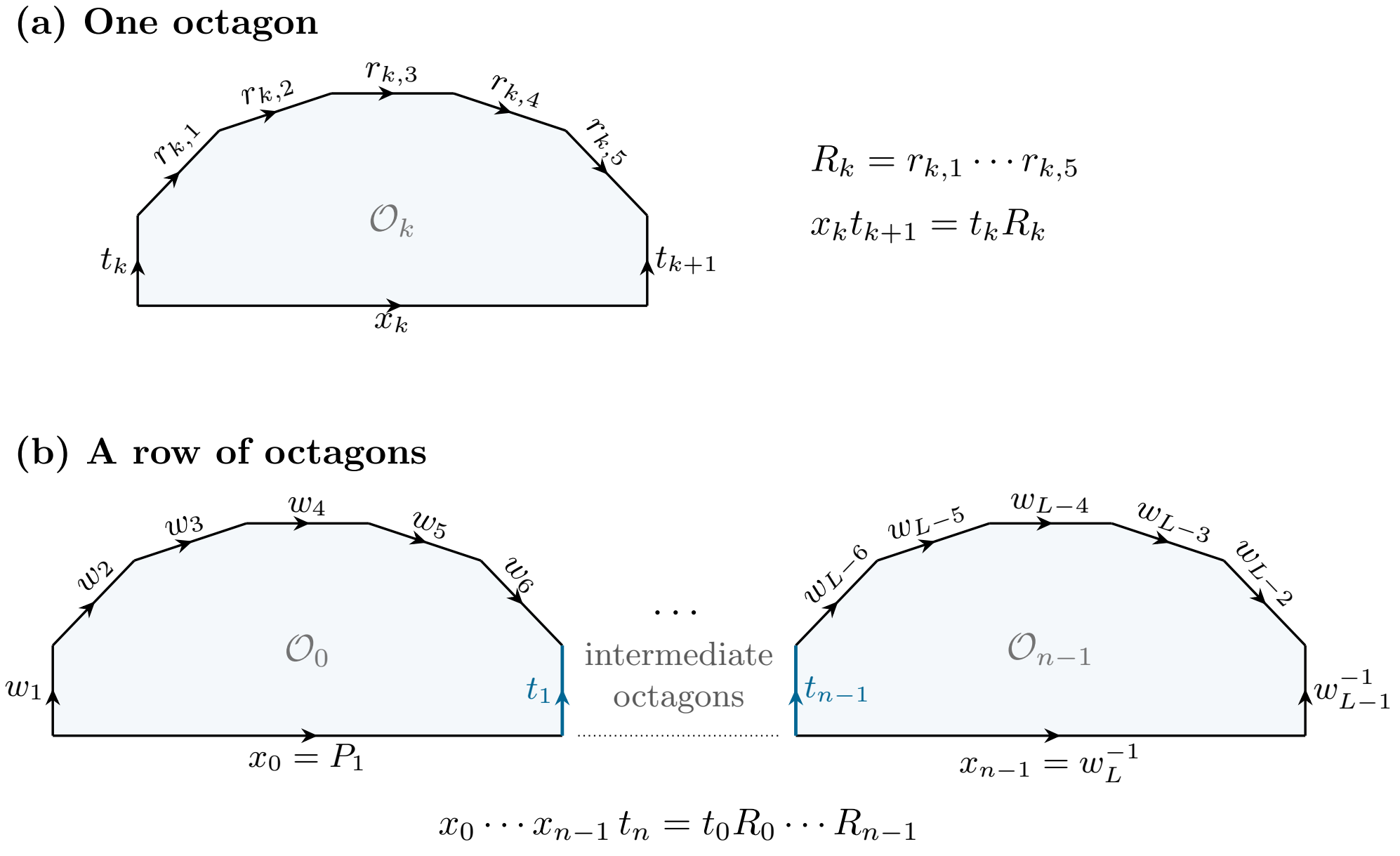
\begin{figure}[htbp]
 \centering
 \begin{tikzpicture}
  \begin{scope}[xshift=1.3cm]
   \relationoctagon{4.8}{x_k}
    {r_{k,1},r_{k,2},r_{k,3},r_{k,4},r_{k,5}}{\mathcal O_k}
   \relationwall{0}{t_k}{left}{}
   \relationwall{4.8}{t_{k+1}}{right}{}
   \node[anchor=west,align=left,font=\small] at (6.2,1.05)
    {$R_k=r_{k,1}\cdots r_{k,5}$\\[5pt]
     $x_k t_{k+1}=t_kR_k$};
  \end{scope}
  \node[anchor=west,font=\small\bfseries] at (0,2.65) {(a) One octagon};
  \begin{scope}[yshift=-4.05cm,xshift=.5cm]
   \relationoctagon{4.8}{x_0=P_1}{w_2,w_3,w_4,w_5,w_6}{\mathcal O_0}
   \relationwall{0}{w_1}{left}{}
   \relationwall{4.8}{t_1}{left}{relation wall}
   \begin{scope}[xshift=7cm]
    \relationoctagon{4.8}{x_{n-1}=w_L^{-1}}
     {w_{L-6},w_{L-5},w_{L-4},w_{L-3},w_{L-2}}{\mathcal O_{n-1}}
    \relationwall{0}{t_{n-1}}{right}{relation wall}
    \relationwall{4.8}{w_{L-1}^{-1}}{right}{}
   \end{scope}
   \draw[densely dotted] (4.95,0)--(6.85,0);
   \node at (5.9,1.15) {$\cdots$};
   \node[font=\footnotesize,align=center,text=black!65] at (5.9,.55)
    {intermediate\\octagons};
   \node[font=\small] at (5.9,-.85)
    {$x_0\cdots x_{n-1}\,t_n=t_0R_0\cdots R_{n-1}$};
  \end{scope}
  \node[anchor=west,font=\small\bfseries] at (0,-1.4)
   {(b) A row of octagons};
 \end{tikzpicture}
 \caption{A generic octagon row, before identifying vertices and labelled edges. Each roof has five edges. Adjacent octagons share an internal wall $t_k$. The outside boundary of the row reads $PW^{-1}$.}
 \label{fig:generic-row}
\end{figure}

We first describe a recipe for gluing octagons given some input words and letters. Let $P=P_1\dots P_p$ be a word, let $n\coloneqq p+1$ and $L\coloneqq 5n+3$, and let $W=w_1\dots w_{L}$ be a positive word. Define
\begin{align*}
    (x_0,\dots,x_{n-1})&\coloneqq(P_1,\dots,P_p,w_L^{-1}),\\
    t_0&\coloneqq w_1,\qquad t_n\coloneqq w_{L-1}^{-1},\\
    r_{k,j}&\coloneqq w_{5k+j+1}\qquad(0\leq k<n,\ 1\leq j\leq 5),
\end{align*}
and let $t_1,\dots,t_{n-1}$ be \emph{local} labels. For each $k=0,\dots, n-1$, consider one octagon $\cO_k$ with cyclic boundary word
\begin{equation*}
    x_kt_{k+1}r_{k,5}^{-1}r_{k,4}^{-1}r_{k,3}^{-1}r_{k,2}^{-1}r_{k,1}^{-1}t_k^{-1},
\end{equation*}
oriented in the counterclockwise direction. Then identify all the vertices together, and edges with the same oriented label. See \Cref{fig:generic-row}, where the row is shown before identifying vertices and edges.

Apply this construction to the following two families of words
\begin{equation}\label{eq:row-inputs}
    \begin{array}{c|c|c|c}
        \text{relation family}&P&W&n\\ \hline
        C_i&c_i b&bC_i&3\\
        D_{ij}&c_i^{-1}d_jc_i&D_{ij}&4.
    \end{array}
\end{equation}
Recall that $\abs{C_i}=17$ and $\abs{D_{ij}}=23$, hence $\abs{W}=5n+3$ in both cases. This yields a collection of rows, each containing three or four octagons. Identify all vertices together and all oriented edges with the same label, preserving their orientation, without identifying the internal labels $t_1,\dots,t_{n-1}$: these are \emph{local} labels in each row, referred to as \emph{internal walls}, and distinct from the internal labels of every other row. All the other edges, labelled $b,c_i,d_j$, are called \emph{old}. Let $Z$ denote the resulting hyperbolic octagon complex.

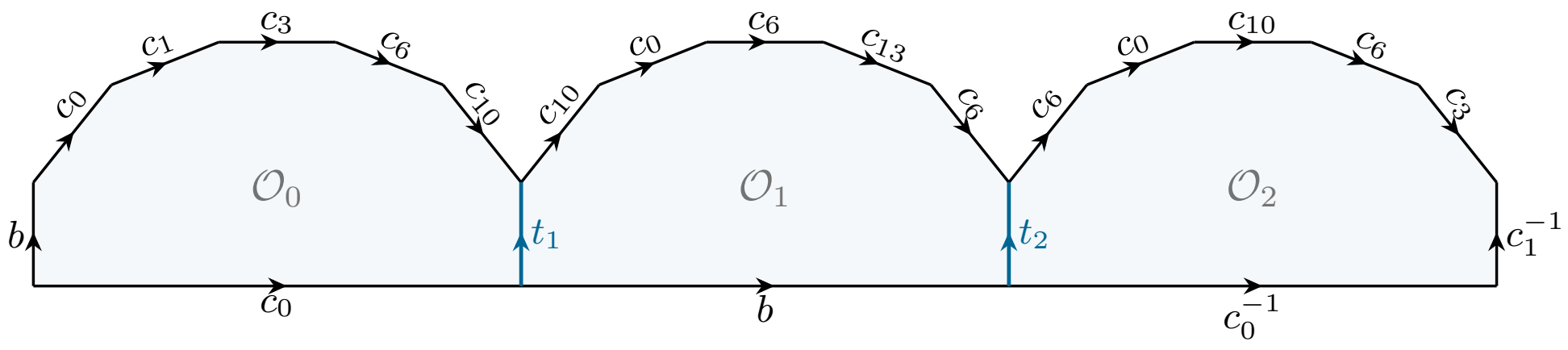
\begin{figure}[htbp]
 \centering
 \begin{tikzpicture}
  \relationoctagon{4}{c_0}{c_0,c_1,c_3,c_6,c_{10}}{\mathcal O_0}
  \begin{scope}[xshift=4cm]
   \relationoctagon{4}{b}{c_{10},c_0,c_6,c_{13},c_6}{\mathcal O_1}
  \end{scope}
  \begin{scope}[xshift=8cm]
   \relationoctagon{4}{c_0^{-1}}{c_6,c_0,c_{10},c_6,c_3}{\mathcal O_2}
  \end{scope}
  \relationwall{0}{b}{left}{}
  \relationwall{4}{t_1}{right}{relation wall}
  \relationwall{8}{t_2}{right}{relation wall}
  \relationwall{12}{c_1^{-1}}{right}{}
 \end{tikzpicture}
 \caption{The $C_i$ family of \eqref{eq:row-inputs}, here shown for $i=0$, hence $P=c_0b$ and $W=bC_0$. Reading the roof from left to right gives the first fifteen letters of $C_0$; its final letters $c_1c_0$ appear on the right wall and last bottom edge with reversed orientations. The outside boundary is $c_0bC_0^{-1}b^{-1}$, yielding the relation $c_0b=bC_0$, or $b^{-1}c_0b=C_0$.  As in \Cref{fig:generic-row}, this is the
 row before the identifications.}
 \label{fig:c-row}
\end{figure}

\begin{figure}[htbp]
 \centering
 \begin{tikzpicture}
  \relationoctagon{3.35}{c_0^{-1}}{d_1,d_3,d_6,d_{10},d_{15}}{\mathcal O_0}
  \begin{scope}[xshift=3.35cm]
   \relationoctagon{3.35}{d_0}{d_{15},d_{21},d_{28},d_{36},d_{45}}{\mathcal O_1}
  \end{scope}
  \begin{scope}[xshift=6.7cm]
   \relationoctagon{3.35}{c_0}{d_{45},d_{55},d_{66},d_{78},d_{91}}{\mathcal O_2}
  \end{scope}
  \begin{scope}[xshift=10.05cm]
   \relationoctagon{3.35}{d_{190}^{-1}}{d_{91},d_{105},d_{120},d_{136},d_{153}}{\mathcal O_3}
  \end{scope}
  \relationwall{0}{d_0}{left}{}
  \relationwall{3.35}{t_1}{right}{relation wall}
  \relationwall{6.7}{t_2}{right}{relation wall}
  \relationwall{10.05}{t_3}{right}{relation wall}
  \relationwall{13.4}{d_{171}^{-1}}{right}{}
 \end{tikzpicture}
 \caption{The $D_{ij}$ family of \eqref{eq:row-inputs}, here shown for $(i,j)=(0,0)$, hence $P=c_0^{-1}d_0c_0$ and $W=D_{00}$. The first letter $d_0$ lies on the left wall, the next twenty letters form the roof, and the final letters $d_{171}d_{190}$ lie on the right wall and last bottom edge with reversed orientations. The boundary word gives the relation $c_0^{-1}d_0c_0=D_{00}$. The internal labels $t_1,t_2,t_3$ are \emph{local} to this row and are distinct from the labels in \Cref{fig:c-row}.}
 \label{fig:d-row}
\end{figure}
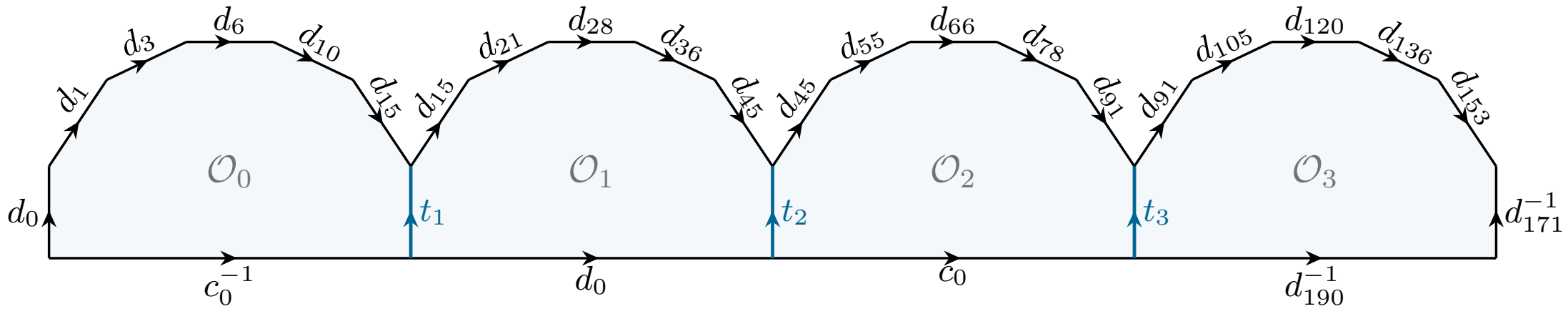

See \Cref{fig:c-row,fig:d-row} for examples of the rows constructed from the relations \eqref{eq:row-inputs}, using the words of \Cref{sec:words_adapted_to_octagons} (in these examples, the words are \eqref{eq:example_C0_D00}). Notice that each pair of repeated letters, coming from zero entries in the tuples $\Delta$ and $\Delta_i$, is separated by an internal wall and thus does not form a corner of an octagon. This will be used in the link description below.

\begin{proposition}\label{prop:combinatorial_description_Z}
    The following assertions hold:
    \begin{enumerate}[label=(\roman*)]
        \item The natural presentation of $\pi_1(Z)$, with old edges as generators and rows of octagons as relations, identifies it with $G$.
        \item Every edge of $Z$ belongs to an octagon and every octagon boundary is immersed in the one-skeleton.
        \item The subgraph of the vertex link spanned by the old half-edges has exactly the following edges:
        \begin{equation}\label{eq:old_link}
            \begin{gathered}
                [b]_+\text{--}[c_i]_+,\qquad
                [b]_-\text{--}[c_i]_+,\qquad
                [c_i]_-\text{--}[c_j]_+\quad(i\ne j),\\
                [c_i]_-\text{--}[d_j]_+,\qquad
                [d_j]_-\text{--}[d_k]_+\quad(j\ne k),
            \end{gathered}
        \end{equation}
        with each edge listed occurring exactly once, and the dash -- denoting an unoriented edge in the vertex link.
        \item The vertex link of $Z$ is simple and triangle-free.
        \item The octagon complex $Z$ is locally \cat{-1}\ and $G$ is torsion-free and hyperbolic.
    \end{enumerate}
\end{proposition}
\begin{proof}
    In the generic recipe for rows of octagons (see \Cref{fig:generic-row}), the relation induced from the boundary of the row is, after eliminating the internal $t_k$'s,
    \begin{equation*}
        x_0\cdots x_{n-1}t_n(R_0\cdots R_{n-1})^{-1}t_0^{-1}=Pw_L^{-1}w_{L-1}^{-1}\cdots w_1^{-1}=PW^{-1}.
    \end{equation*}
    In the present construction, with \eqref{eq:row-inputs} as inputs, the resulting relations are $c_ib=bC_i$ and $c_i^{-1}d_jc_i=D_{ij}$. This corresponds exactly to the presentation of $G$, hence the natural map $G\to\pi_1(Z)$ sending letters to the corresponding old edges is an isomorphism. This proves (i).
    
    For (ii), the only assertion to verify is that boundaries of octagons are immersed in the one-skeleton of $Z$. This follows from the fact that the boundary words of octagons have no consecutive letters with the same unoriented label: consecutive identical letters in the $C_i$ and $D_{ij}$ are separated by internal walls, so they do not form a corner of any octagon. The same holds for the bottom and wall corners by inspection of \Cref{fig:c-row,fig:d-row}.

    Let $\Link(v)$ be the vertex link of the single vertex $v$ of $Z$. Edges of $\Link(v)$ correspond bijectively to corners of octagons in $Z$, of which there are two types: corners between old half-edges (\emph{old corners}), and corners involving an internal wall (\emph{internal corners}).

    The subgraph of old corners, or \emph{old graph}, can be described using \Cref{lem:words_injection}. Old corners in the roofs and at the bottom-right of the rows correspond precisely to ordered pairs of distinct $c$-letters and of distinct $d$-letters, with a pair $c_ic_j$ contributing exactly the edge $[c_i]_{-}\text{--}[c_j]_{+}$, and similarly for a pair of distinct $d$-letters. The rest of the old corners (the ones involving the leftmost wall) in the $C$-rows and $D$-rows contribute to the other edges listed in \eqref{eq:old_link}, each appearing exactly once by \Cref{lem:words_injection} and its proof. In particular, this uses the fact that $c_i$ is the starting letter of the word $C_i$, and that $d_j$ is the starting letter of each of the $D_{ij}$'s. This yields \eqref{eq:old_link} as the exhaustive and non-redundant list of edges in the subgraph on old half-edges, proving (iii).

    There is no corner between two edges with the same label, hence $\Link(v)$ is loop-free. There are no parallel edges in the old graph by (iii), and none involving an internal wall, since each $t_k$ lies in exactly four corners with four distinct neighboring half-edges. Thus $\Link(v)$ is simple. It remains to show that $\Link(v)$ is triangle-free.

    From \eqref{eq:old_link}, it follows that the old graph is bipartite with parts
    \begin{equation}\label{eq:actual_bipartition}
        \{[c_i]_+,[d_j]_+\},\qquad \{[b]_+,[b]_-,[c_i]_-,[d_j]_-\}.
    \end{equation}
    Thus, the old graph is triangle-free. A triangle cannot involve $[t_k]_-$ for an internal wall $t_k$, as otherwise there would be an edge $[c_i]_-\text{--}[c_i]_+$ or $[d_j]_-\text{--}[d_j]_+$ (see \Cref{fig:c-row,fig:d-row}), not possible by (iii). Similarly, a triangle involving $[t_k]_+$ for an internal wall $t_k$ of a $C$-row would require the existence of an edge $[c_i]_-\text{--}[b]_+$ or $[b]_-\text{--}[c_i]_-$. Finally, a triangle involving $[t_k]_+$ for an internal wall $t_k$ of a $D$-row would require the existence of an edge $[c_i]_+\text{--}[d_j]_+$, $[d_j]_-\text{--}[c_i]_+$ or $[c_i]_-\text{--}[d_{j'}]_-$ where $d_{j'}$ is the last letter of $D_{ij}$. None of these link edges exist by (iii) again. This shows that $\Link(v)$ is triangle-free, proving (iv).

    Finally, (v) follows from (iv) by \Cref{prop:hyperbolicity_of_octagon_complex}.
\end{proof}

\subsection{Attachment criterion}\label{sec:attachment_criterion}

The next step is to attach a genus-$431$ surface to $Z$ along a bouquet of circles representing $F$. We require $F$ to be quasi-convex in the surface group, so that $\partial F$ embeds in the boundary of the surface group and the convergence $h_nd_F^+\to b_F^+$ persists there. At the same time, the attachment must preserve injectivity of the induced homomorphisms between fundamental groups and the link conditions needed to apply the theorem of Desgroseilliers--Haglund (\Cref{thm:Desgroseilliers--Haglund}). The following criterion ensures these properties; see also~\cite[Proof of Theorem 3.1]{MR2308946}.

\begin{lemma}[Attachment criterion]\label{lem:attachment_criterion}
    Let $Z'$ and $Y'$ be finite connected complexes made from the same regular right-angled octagon. Suppose their vertex links are simple and triangle-free, that every edge belongs to a face and that all face boundaries are immersed. Let $R'$ be a bouquet of circles, with embeddings $f\colon R'\monic Z'$ and $g\colon R'\monic Y'$ whose images lie in the one-skeleta, mapping edges to edges and injective on half-edges. Call the half-edges in the image of these embeddings \emph{marked half-edges}. Assume that $f_*\colon\pi_1(R')\to\pi_1(Z')$ is injective and that the following hold in the vertex links in $Y'$:
    \begin{enumerate}[label=(\roman*)]
        \item distinct marked half-edges are at distance at least two in vertex links in $Y'$;
        \item if a link path made of two edges joins marked half-edges in $Y'$, then the corresponding marked half-edges are not adjacent in the link of the corresponding vertex in $Z'$.
    \end{enumerate}
    Then the attachment $Q'\coloneqq Z'\cup_{R'}Y'$ has simple triangle-free links, immersed face boundaries and every edge in a face. It is locally \cat{-1}. With $F'\coloneqq\pi_1(R')$, one has
    \begin{equation}\label{eq:attachment_criterion_SvK}
        \pi_1(Q')=\pi_1(Z')\ast_{F'}\pi_1(Y'),
    \end{equation}
    with both factors embedded in $\pi_1(Q')$. Moreover, the maps $g\colon R'\to Y'$ and $Z'\to Q'$ lift to convex isometric embeddings of universal covers, $F'$ is quasi-convex in $\pi_1(Y')$ and $\pi_1(Z')$ is quasi-convex in $\pi_1(Q')$.
\end{lemma}
\begin{proof}
At each vertex of $R'$, the vertex link in $Q'$ is the union of the links in $Z'$ and $Y'$ along the marked half-edges. By (i), no edge is shared and no new triangle has exactly one $Y'$-edge. Thus a new triangle would have two $Y'$-edges and one $Z'$-edge, which is excluded by (ii). Hence links in $Q'$ are simple and triangle-free, and $Q'$, like $Y'$ and $Z'$, is locally \cat{-1}\ by \Cref{prop:hyperbolicity_of_octagon_complex}. Face boundaries remain immersed and every edge still lies in a face.

All link edges have length $\pi/2$. Under this condition, the inclusion of a subcomplex $A\subseteq K$ is a local isometry at a vertex $v$ if and only if $\Link_A(v)$ is a full subgraph of $\Link_K(v)$, i.e.\ if vertices $e,e'$ in $\Link_A(v)$ form an edge in $\Link_K(v)$, then they form an edge in $\Link_A(v)$ (see~\cite[Lemma 2.6]{MR3110583}, or~\cite[\S 2.2]{MR2377497} for the same statement for cube complexes).

By (i), no two marked half-edges are adjacent in the links of $Y'$. Thus the link maps induced by $g$ have full image. Similarly, the links of $Z'$ are full subgraphs of the links of $Q'$: any additional link edge would come from $Y'$ with both half-edges marked, contradicting (i). Local convexity also holds at interior points of edges. Thus, $g$ and $Z'\to Q'$ are local isometric embeddings.

By \Cref{prop:local_isometries_are_pi_1_injective}, both maps lift to isometric embeddings of universal covers and $g_*$ and $\pi_1(Z')\to\pi_1(Q')$ are injective. Both lifts have convex image: their targets are \cat{-1}, hence uniquely geodesic by~\cite[II.1.4(1)]{MR1744486}, and an isometric embedding of a geodesic space into a uniquely geodesic space has convex image. Since $f_*$ is assumed to be injective and $g_*$ is injective by the above, Seifert--van Kampen and the normal form \cite[IV, \S 2]{MR577064} yield \eqref{eq:attachment_criterion_SvK} with both factors embedded. Cocompactness on the convex lifts yields quasi-convexity.
\end{proof}

\subsection{Gluing of the genus-431 surface}\label{sec:gluing_the_surface}

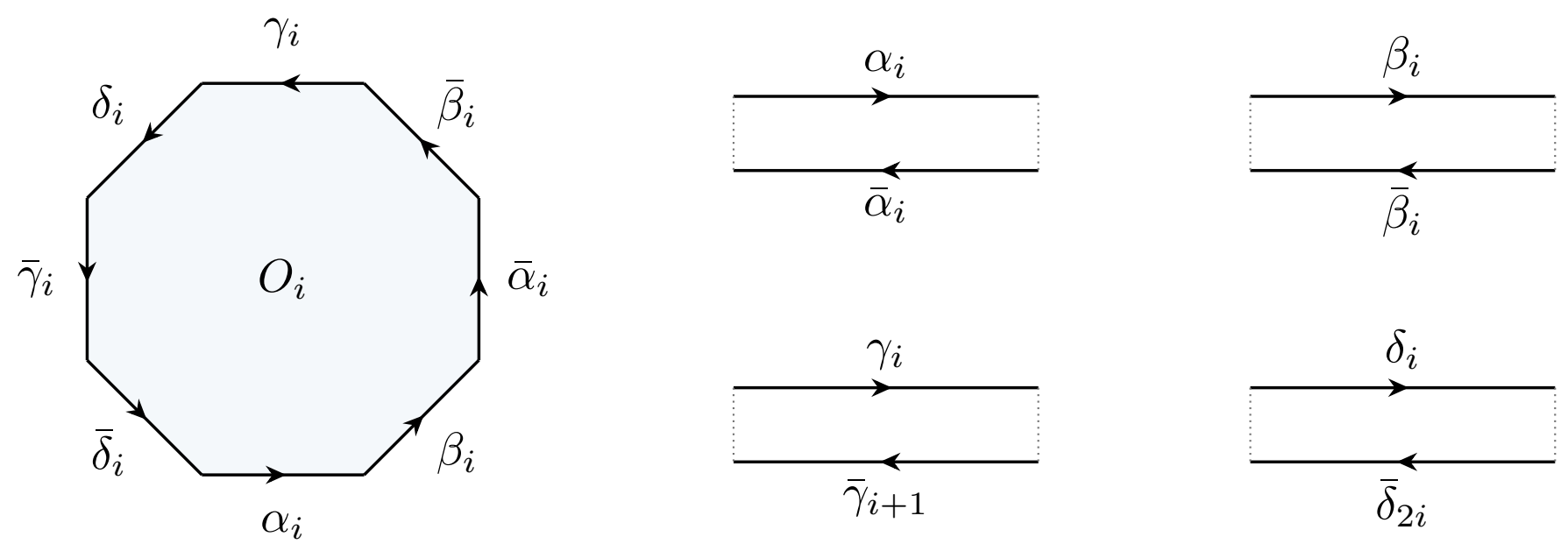
\begin{figure}[htbp]
 \centering
\begin{tikzpicture}[>=Stealth]
  \begin{scope}[xshift=2.1cm]
   \occurrenceoctagon{1.6}
    {\alpha_i/pairalpha,\beta_i/pairbeta,\bar\alpha_i/pairalpha,\bar\beta_i/pairbeta,
     \gamma_i/pairgamma,\delta_i/pairdelta,\bar\gamma_i/pairgamma,\bar\delta_i/pairdelta}
    {O_i}
  \end{scope}
  \begin{scope}[xshift=5.5cm,yshift=1.1cm]
   \sidepair{\alpha_i}{\bar\alpha_i}{pairalpha}
  \end{scope}
  \begin{scope}[xshift=9.4cm,yshift=1.1cm]
   \sidepair{\beta_i}{\bar\beta_i}{pairbeta}
  \end{scope}
  \begin{scope}[xshift=5.5cm,yshift=-1.1cm]
   \sidepair{\gamma_i}{\bar\gamma_{i+1}}{pairgamma}
  \end{scope}
  \begin{scope}[xshift=9.4cm,yshift=-1.1cm]
   \sidepair{\delta_i}{\bar\delta_{2i}}{pairdelta}
  \end{scope}
 \end{tikzpicture}
 \caption{One octagon $O_i$, and the gluing pattern between the octagons, with $i\in\Z/T\Z$.}
 \label{fig:surface-gluing}
\end{figure}

Let $T=287$. For each $i\in\Z/T\Z$, consider an oriented regular right-angled octagon with boundary sides
\begin{equation*}
    \alpha_i,\ \beta_i,\ \bar\alpha_i,\ \bar\beta_i,\ \gamma_i,\ \delta_i,\ \bar\gamma_i,\ \bar\delta_i
\end{equation*}
in the counterclockwise order. Glue the sides of these octagons according to the pairing (see \Cref{fig:surface-gluing})
\begin{equation}\label{eq:surface_pairing}
    \alpha_i\leftrightarrow\bar\alpha_i,\qquad \beta_i\leftrightarrow\bar\beta_i,\qquad \gamma_i\leftrightarrow\bar\gamma_{i+1},\qquad \delta_i\leftrightarrow\bar\delta_{2i},
\end{equation}
each side being identified with another side by an orientation-reversing isometry.

Note that $i\mapsto 2i$ is a bijection since $T$ is odd. Denote the quotient of this gluing by $Y_0$. Orient each edge by the orientation of the unbarred side before identification.

Let $R$ be the subcomplex of the one-skeleton of $Y_0$ made of the $\delta$-edges; it is an embedded bouquet of 287 circles. Consider the embedding $m\colon R\monic Z$ defined by
\begin{equation}\label{eq:surface_marking}
    \delta_0\longmapsto b^{-1},\qquad\delta_i\longmapsto d_{i-1}^{\,(-1)^{i+1}}\quad(1\le i<T).
\end{equation}

\begin{proposition}\label{prop:construction_of_Q_0}
    The quotient $Y_0$ is a closed connected oriented surface of genus 431, with one vertex, 1148 edges and 287 octagons. Its $\delta$-edges form an embedded bouquet $R$ in $Y_0$, which satisfies the hypotheses of \Cref{lem:attachment_criterion} with inclusion $R\monic Y_0$ and embedding $m\colon R\monic Z$ defined by \eqref{eq:surface_marking}. Thus, letting
    \begin{equation*}
        Q_0=Z\cup_RY_0,\qquad \Gamma_{431}=\pi_1(Y_0),\qquad K_0=\pi_1(Q_0),
    \end{equation*}
    we have $F= m_*(\pi_1(R))$, $K_0=G\ast_{F}\Gamma_{431}$ and both $G=\pi_1(Z)$ and $\Gamma_{431}$ embed in $K_0$. Via the identification $F\cong\pi_1(R)$, $F$ is identified with a quasi-convex subgroup of $\Gamma_{431}$.

    Moreover, $Q_0$ is a finite one-vertex locally \cat{-1}\ octagon complex with immersed face boundaries, every edge in a face, and simple triangle-free vertex link. Its fundamental group $K_0$ is torsion-free and hyperbolic.
\end{proposition}
\begin{proof}
    Let us describe the link of a vertex $v$ in $Y_0$, directly from the edge identifications \eqref{eq:surface_pairing} (see \Cref{fig:surface-gluing}), starting at the half-edge $[\alpha_i]_+$ (recall that $+$ denotes outgoing, or positive, half-edges):
    \begin{equation}\label{eq:link_vertex_Y0}
        \dots\text{--}\underbrace{[\alpha_i]_+\text{--}[\bar\beta_{i}]_+\text{--}[\bar\alpha_{i}]_+\text{--}[\beta_{i}]_+\text{--}[\gamma_{i}]_+\text{--}[\bar\delta_{i+1}]_+\text{--}[\bar\gamma_{\tfrac{i+1}{2}}]_+\text{--}[\delta_{\tfrac{i-1}{2}}]_+}_{\text{repeating block}}\text{--}[\alpha_{i-1}]_+\text{--}\dots,
    \end{equation}
    where the repeating block of size $8$ is repeated with a shift of index $i\mapsto i-1$, until it closes; see \Cref{fig:surface_vertex_neighborhood}. Division by $2$ is well-defined in $\Z/T\Z$. Observe that \eqref{eq:link_vertex_Y0} contains all positive half-edges, and each negative half-edge is identified to a positive half-edge under \eqref{eq:surface_pairing}. Thus, $v$ is the only vertex of $Y_0$ and $\Link(v)$ is a single cycle of length $8T=2296$. In particular, $Y_0$ is a closed oriented surface, endowed with a hyperbolic cone metric, with a single cone point of angle $1148\pi\geq 2\pi$. Hence, $Y_0$ is locally \cat{-1}; this also follows from \Cref{prop:hyperbolicity_of_octagon_complex}. Note also that $Y_0$ has immersed face boundaries, and each edge lies in the boundary of a face by construction.

    \begin{figure}[htbp]
    \centering
        \begin{tikzpicture}[scale=.9,transform shape,line cap=round,
        corner side/.style={line width=.65pt,postaction={decorate},
        decoration={markings,mark=at position .6 with
            {\arrow{Stealth[length=4pt,width=3.6pt]}}}},
        germ ray/.style={MidnightBlue,line width=.7pt,postaction={decorate},
        decoration={markings,mark=at position .64 with
            {\arrow{Stealth[length=4pt,width=3.6pt]}}}}]
        \node[font=\small\bfseries] at (0,3.1) {(a) Before identification};
        \foreach \k/\incoming/\outgoing in {
        -1/{\bar\delta_i}/{\alpha_i},
        0/{\bar\alpha_i}/{\bar\beta_i},
        1/{\beta_i}/{\bar\alpha_i},
        2/{\alpha_i}/{\beta_i},
        3/{\bar\beta_i}/{\gamma_i},
        4/{\bar\gamma_{i+1}}/{\bar\delta_{i+1}},
        5/{\delta_{\frac{i+1}{2}}}/{\bar\gamma_{\frac{i+1}{2}}},
        6/{\gamma_{\frac{i-1}{2}}}/{\delta_{\frac{i-1}{2}}},
        7/{\bar\delta_{i-1}}/{\alpha_{i-1}}} {
        \pgfmathsetmacro{\midangle}{195-30*\k}
        \pgfmathsetmacro{\startangle}{\midangle+13}
        \pgfmathsetmacro{\endangle}{\midangle-13}
        \coordinate (c\k) at (\midangle:.32);
        \begin{scope}[shift={(c\k)}]
            \fill[black!5] (0,0)--(\startangle:1.58)
            arc[start angle=\startangle,end angle=\endangle,radius=1.58]--cycle;
            \draw[black!45,dashed,line width=.4pt] (\startangle:.59)
            arc[start angle=\startangle,end angle=\endangle,radius=.59];
            \draw[corner side,black!75] (\startangle:1.71)--(0,0);
            \draw[corner side,MidnightBlue] (0,0)--(\endangle:1.71);
            \pgfmathsetmacro{\incominglabelangle}{\startangle-3}
            \pgfmathsetmacro{\incomingrotation}{\startangle>90 ? \startangle-180 : \startangle}
            \pgfmathtruncatemacro{\incomingflip}{\startangle>90 ? 1 : 0}
            \ifnum\incomingflip=1
            \node[rotate=\incomingrotation,anchor=east,font=\footnotesize,inner sep=1.5pt]
            at (\incominglabelangle:1.82) {$\incoming$};
            \else
            \node[rotate=\incomingrotation,anchor=west,font=\footnotesize,inner sep=1.5pt]
            at (\incominglabelangle:1.82) {$\incoming$};
            \fi
            \pgfmathsetmacro{\outgoinglabelangle}{\endangle+3}
            \pgfmathsetmacro{\outgoingrotation}{\endangle>90 ? \endangle-180 : \endangle}
            \pgfmathtruncatemacro{\outgoingflip}{\endangle>90 ? 1 : 0}
            \ifnum\outgoingflip=1
            \node[rotate=\outgoingrotation,anchor=east,font=\footnotesize,inner sep=1.5pt,
            text=MidnightBlue] at (\outgoinglabelangle:1.82) {$\outgoing$};
            \else
            \node[rotate=\outgoingrotation,anchor=west,font=\footnotesize,inner sep=1.5pt,
            text=MidnightBlue] at (\outgoinglabelangle:1.82) {$\outgoing$};
            \fi
            \fill (0,0) circle[radius=1pt];
        \end{scope}
        }
        \node[font=\Large] at (0,-1.43) {$\cdots$};
        \begin{scope}[xshift=6.8cm]
            \node[font=\small\bfseries] at (0,3.1) {(b) After identification};
            \fill[black!2] (0,0) circle[radius=1.9];
            \foreach \k in {0,...,7} {
                \pgfmathsetmacro{\startangle}{210-30*\k}
                \pgfmathsetmacro{\endangle}{\startangle-30}
                \pgfmathtruncatemacro{\tint}{mod(\k,2)==0 ? 7 : 3}
                \fill[MidnightBlue!\tint] (0,0)--(\startangle:1.9)
                arc[start angle=\startangle,end angle=\endangle,radius=1.9]--cycle;
                \draw[black!35,dashed,line width=.3pt] (\startangle:.43)
                arc[start angle=\startangle,end angle=\endangle,radius=.43];
            }
            \draw[black!45,dashed,line width=.4pt] (0,0) circle[radius=1.9];
            \foreach \ang/\germ/\anch in {
                210/{[\alpha_i]_+}/north east,
                180/{[\bar\beta_i]_+}/east,
                150/{[\bar\alpha_i]_+}/south east,
                120/{[\beta_i]_+}/south,
                90/{[\gamma_i]_+}/south,
                60/{[\bar\delta_{i+1}]_+}/south west,
                30/{[\bar\gamma_{\frac{i+1}{2}}]_+}/south west,
                0/{[\delta_{\frac{i-1}{2}}]_+}/west,
                -30/{[\alpha_{i-1}]_+}/north west} {
                \draw[germ ray] (0,0)--(\ang:2.02);
                \node[font=\small,text=MidnightBlue,anchor=\anch,inner sep=1.5pt]
                at (\ang:2.13) {$\germ$};
            }
            \fill (0,0) circle[radius=1.4pt];
            \node[font=\small,anchor=north,inner sep=3pt] at (0,0) {$v$};
            \node[font=\Large] at (0,-1.45) {$\cdots$};
        \end{scope}
        \end{tikzpicture}
    \caption{Consecutive corners before (a) and after (b) identification
    at the unique vertex of $Y_0$. In (a), facing parallel sides are paired under \eqref{eq:surface_pairing}. In (b), the sides are paired and the positive half-edges appear in the order \eqref{eq:link_vertex_Y0}. Each corner has angle $\pi/2$ in $Y_0$.}
    \label{fig:surface_vertex_neighborhood}
    \end{figure}
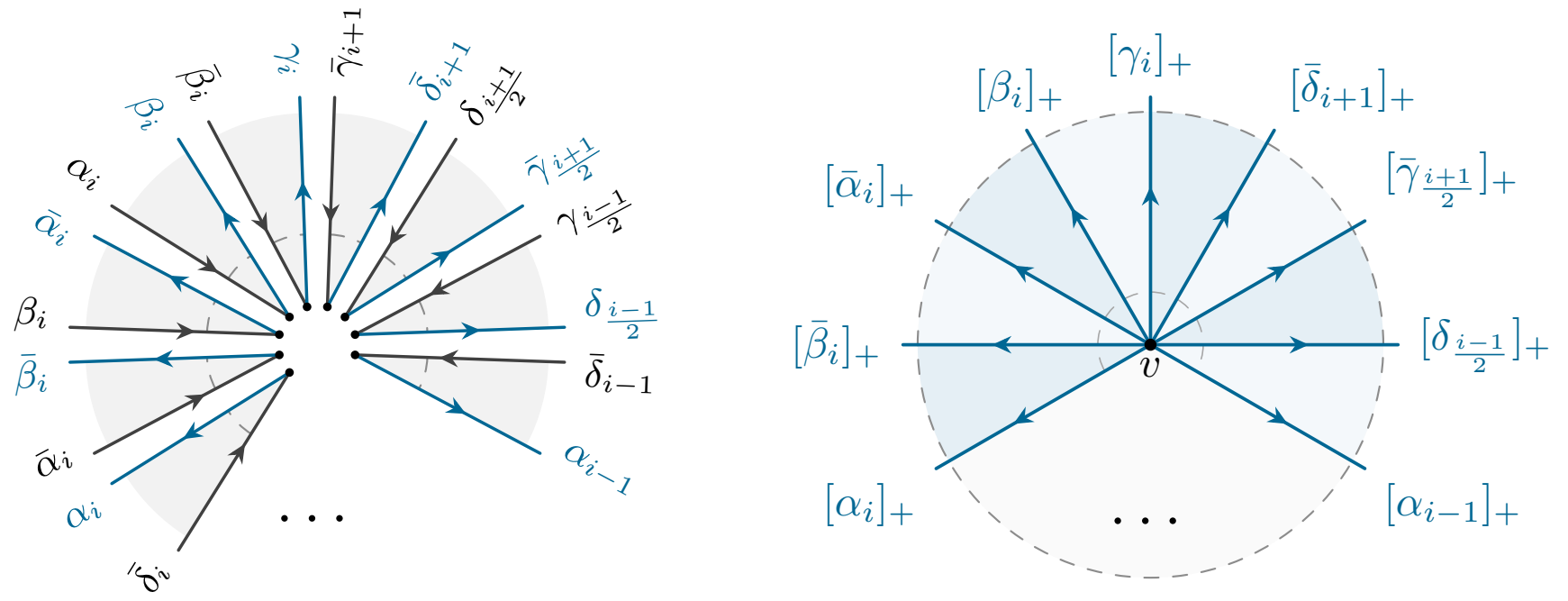
    
    The complex $Y_0$ has $1$ vertex, $4T=1148$ edges ($8T$ before the identification) and $T=287$ octagonal faces. Let $g$ be its genus. An Euler characteristic computation yields
    \begin{equation*}
        2-2g=1-4T+T=-860,\qquad\text{hence }\quad g=431.
    \end{equation*}

    It remains to show that the bouquet $R$ made of the $\delta$-edges and its inclusions $R\monic Y_0$ and $m\colon R\monic Z$ satisfy the conditions of the attachment criterion (\Cref{lem:attachment_criterion}): that $m_*\colon\pi_1(R)\to\pi_1(Z)$ is injective and the two conditions on the links. Injectivity follows from \Cref{lem:double_HNN}, which says that $\langle b,D\rangle=m_*(\pi_1(R))$ is free in $G=\pi_1(Z)$.

    From the vertex link description \eqref{eq:link_vertex_Y0}, see also \Cref{fig:surface_vertex_neighborhood}, $\delta$-half-edges in $\Link(v)$ are separated by $2$ or $6$ link edges, i.e.\ two or six octagon corners, thus condition (i) of \Cref{lem:attachment_criterion} is satisfied. For (ii), the only pairs of marked half-edges at link distance 2 in $Y_0$ are the endpoints of the two-edge link path
    \begin{equation*}
        [\bar\delta_{i+1}]_+\text{--}[\bar\gamma_{\tfrac{i+1}{2}}]_+\text{--}[\delta_{\tfrac{i-1}{2}}]_+,\qquad\text{ for }i\in\Z/T\Z,
    \end{equation*}
    which, after the identification \eqref{eq:surface_pairing} and the reindexing $i\mapsto (i-1)/2$, are
    \begin{equation}\label{eq:check_attachment_potential_pair}
        [\delta_{i+1}]_-,\quad [\delta_{i}]_+,\qquad\text{ for }i\in\Z/T\Z.
    \end{equation}
    First, assume that $1\leq i\leq T-2$. Then this pair is sent via $m$ (see \eqref{eq:surface_marking}) to the pair $\{[d_i^{\mp 1}]_-,[d_{i-1}^{\pm 1}]_+\}$ where $\pm 1$ and $\mp 1$ depend on the parity of $i$ and, importantly, are opposite to each other. This corresponds to the pair $\{[d_i^{\pm 1}]_+,[d_{i-1}^{\pm 1}]_+\}$, which is not a link edge in $Z$ by \Cref{prop:combinatorial_description_Z}(iii). For $i=0$, \eqref{eq:check_attachment_potential_pair} is mapped via $m$ (see \eqref{eq:surface_marking}) to the pair $\{[b]_-,[d_0]_{-}\}$, again not a link edge in $Z$ by \Cref{prop:combinatorial_description_Z}(iii). Finally, for $i=T-1$, the pair maps to $\{[b]_+,[d_{285}]_{-}\}$, also absent from the link of $Z$. Thus, condition (ii) of \Cref{lem:attachment_criterion} is satisfied.

    Therefore, \Cref{lem:attachment_criterion} applies and, by \Cref{prop:hyperbolicity_of_octagon_complex}, $K_0=\pi_1(Q_0)$ is torsion-free and hyperbolic.
\end{proof}

\section{Action on real-hyperbolic space}\label{sec:action_on_real-hyperbolic_space}

The key ingredient of the construction is the following theorem of Desgroseilliers and Haglund, which provides a criterion for the fundamental group of a polygonal complex to act convex cocompactly on a real-hyperbolic space of some (unprescribed) dimension.

\begin{theorem}[Desgroseilliers--Haglund {\cite[Theorem 1.7]{MR3110583}}]\label{thm:Desgroseilliers--Haglund}
    Let $X$ be a simply connected locally compact regular Euclidean $2m$-gonal complex, where $m\geq 4$, with simple vertex links of girth at least four.
    
    Every uniform lattice $K$ of $X$ admits a discrete faithful convex cocompact representation in $\Isom(\HH^q)$ for some positive integer $q$.
\end{theorem}

The above theorem concerns Euclidean polygonal complexes, i.e.\ made of \emph{Euclidean polygons}, so $Q_0$ will be given a different metric in the proof. Also, attaching maps must be local isometries. This condition holds for $Q_0$ since its face boundaries are immersed.

We can now conclude the proof of the main result, whose statement is recalled.

\maintheorem*

\begin{proof}[Proof of \Cref{thm:main}]
    Let $Q_0=Z\cup_R Y_0$ be the octagon complex constructed in \Cref{prop:construction_of_Q_0}, with fundamental group $K_0=G\ast_F\Gamma_{431}$. By \Cref{prop:construction_of_Q_0}, $Q_0$ is a finite octagon complex with immersed face boundaries and whose single vertex has a simple triangle-free link. Also, $K_0$ is hyperbolic and torsion-free.

    To apply \Cref{thm:Desgroseilliers--Haglund}, we endow $Q_0$ with the metric where each of its faces is a regular Euclidean octagon of side length $1$, and let $X\coloneqq\widetilde{Q}_0$ be its universal cover. It is simply connected, locally compact, and its links are simple and triangle-free, hence their combinatorial girth, i.e.\ length of their shortest embedded cycles, is at least $4$.

    The group $K_0$ acts cocompactly on $X$ by automorphisms of the polygonal complex, and the action is properly discontinuous. Thus, $K_0$ is a uniform lattice of $X$.

    By \Cref{thm:Desgroseilliers--Haglund}, there exists $q_0\in\N$ and a discrete faithful convex cocompact representation
    \begin{equation*}
        \sigma_0\colon K_0\to\Isom(\HH^{q_0}),
    \end{equation*}
    which must be purely loxodromic: no elliptic elements since $\sigma_0$ is discrete faithful and $K_0$ is torsion-free, and no parabolic elements since $\sigma_0$ is convex cocompact.

    By \Cref{prop:construction_of_Q_0}, the natural homomorphism $\iota\colon\Gamma_{431}\to K_0$ is injective. Thus the representation $\rho\coloneqq\sigma_0\circ\iota$ is discrete, faithful and purely loxodromic. Since $\sigma_0$ is convex cocompact, there is a strictly positive lower bound on translation lengths of non-trivial elements under $\sigma_0$, hence under $\rho$ too.

    Similarly, $\sigma_0\restr{G}$ is discrete, faithful and purely loxodromic, hence $\sigma_0\restr{F}$ has no \ct\ map by \Cref{prop:explicit_obstruction}, with obstruction given precisely by \eqref{eq:intro_witness}. Since the inclusion $F\monic\Gamma_{431}$ is quasi-convex (\Cref{prop:construction_of_Q_0}) and by \Cref{lem:quasi-convex_quasi-isometric}(ii), it extends to an $F$-equivariant embedding $\partial F\monic\partial\Gamma_{431}$; \eqref{eq:intro_witness} holds for $\rho$ too, thus $\rho$ admits no \ct\ map either. In particular, it cannot be convex cocompact.
\end{proof}

The construction also provides a discrete faithful purely loxodromic action of a free group with no \ct\ map.

\begin{corollary}\label{cor:free_group_action_with_no_CT}
    Restricting $\sigma_0$ from \Cref{thm:main} to $F=\langle b,D\rangle\cong\FF_{287}$ yields a discrete faithful purely loxodromic representation $\FF_{287}\to\Isom(\HH^{q_0})$, for the same $q_0\in\N$, with no \ct\ map.
\end{corollary}

\section{Extension to larger genera}\label{sec:extension_to_large_genera}

From \Cref{thm:main}, one can produce counterexamples in arbitrarily large genera by taking finite covers. However, the construction can be adapted directly to yield counterexamples in genus $g$, for each $g\geq 431$, acting on a real-hyperbolic space of some fixed dimension $q_*$, possibly larger than $q_0$.

\maintheoremhighergenera*

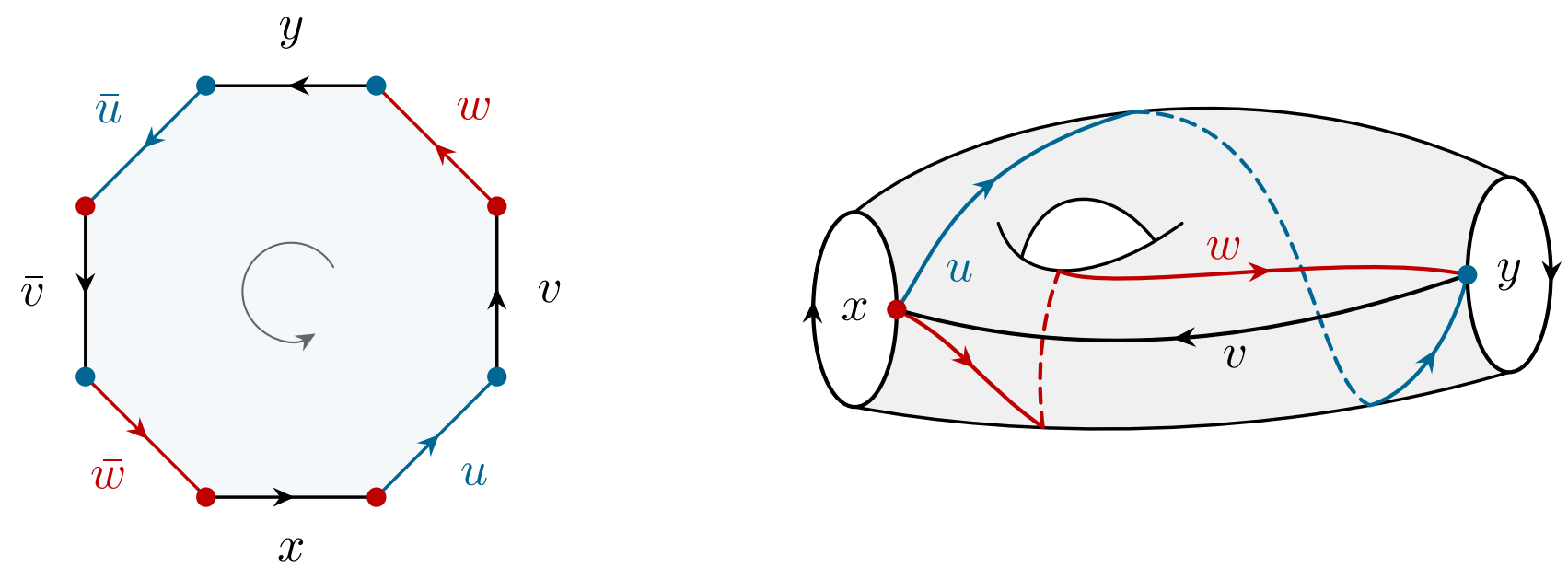
\begin{figure}[htbp]
 \centering
 \begin{tikzpicture}[line cap=round,line join=round,
  handle edge/.style={line width=.8pt,postaction={decorate},
   decoration={markings,mark=at position .52 with
    {\arrow{Stealth[length=4.6pt,width=4pt]}}}},
  handle hidden/.style={densely dashed,line width=.8pt},
  handle label/.style={font=\small,inner sep=1.5pt}]
  \begin{scope}[shift={(1.8,0)}]
   \occurrenceoctagon{1.6}
    {x/black,u/{MidnightBlue},v/black,w/{red!75!black},
     y/black,\bar u/{MidnightBlue},\bar v/black,\bar w/{red!75!black}}{}
   \foreach \k in {0,1,3,6} {\fill[red!75!black] (o\k) circle[radius=2pt];}
   \foreach \k in {2,4,5,7} {\fill[MidnightBlue] (o\k) circle[radius=2pt];}
   \draw[-{Stealth[length=4pt,width=3.6pt]},black!60] (30:.35)
    arc[start angle=30,end angle=300,radius=.35];
  \end{scope}
  \begin{scope}[shift={(8.2,0)}]
   \filldraw[fill=black!6,draw=black,line width=.65pt]
    (-2.35,.57) .. controls (-1.5,1.28) and (.6,1.7) .. (2.35,.82)
    arc[start angle=90,end angle=-90,x radius=.30,y radius=.7]
    .. controls (.6,-1.1) and (-1.2,-1.05) .. (-2.35,-.83)
    arc[start angle=-90,end angle=-270,x radius=.30,y radius=.7]--cycle;
   \coordinate (hp) at (-2.05,-.13);
   \coordinate (hq) at (2.05,.12);
   \draw[handle edge,MidnightBlue] (hp)
    .. controls (-1.8,.20) and (-1.65,.92) .. (-.35,1.29);
   \draw[handle hidden,MidnightBlue] (-.35,1.29)
    .. controls (.9,1.3) and (0.9,-.63) .. (1.35,-.818);
   \draw[handle edge,MidnightBlue] (1.35,-.818)
    .. controls (1.65,-.73) and (1.94,-.35) .. (hq);
   \node[handle label,text=MidnightBlue] at (-1.60,.15) {$u$};
   \draw[handle edge,red!75!black] (hp)
    .. controls (-1.6,-.36) and (-1.45,-.66) .. (-1.0,-.975);
   \draw[handle hidden,red!75!black] (-1.0,-.975)
    .. controls (-1.07,-.40) and (-.95,.08) .. (-.75,.44);
   \fill[white] (-1.15,.25)
    .. controls (-1.00,.80) and (-.50,.76) .. (-.20,.37)
    .. controls (-.50,.14) and (-.90,.06) .. cycle;
   \draw[line width=.65pt] (-1.15,.25)
    .. controls (-1.00,.80) and (-.50,.76) .. (-.20,.37);
   \draw[line width=.65pt] (-1.32,.49)
    .. controls (-1.15,-.04) and (-.50,.12) .. (0.0,.49);
   \draw[handle edge,red!75!black] (-.8837,.1478)
    .. controls (-.52,-.03) and (1.19,.30) .. (hq);
   \node[handle label,text=red!75!black] at (.3,.31) {$w$};
   \draw[handle edge] (hq)
    .. controls (.65,-.37) and (-.70,-.52) .. (hp);
   \node[handle label] at (.38,-.47) {$v$};
   \filldraw[fill=white,draw=black,handle edge] (hp)
    arc[start angle=0,end angle=-360,x radius=.30,y radius=.7]--cycle;
   \filldraw[fill=white,draw=black,handle edge] (hq)
    arc[start angle=180,end angle=-180,x radius=.30,y radius=.7]--cycle;
   \node[handle label] at (-2.35,-.13) {$x$};
   \node[handle label] at (2.35,.12) {$y$};
   \fill[red!75!black] (hp) circle[radius=2pt];
   \fill[MidnightBlue] (hq) circle[radius=2pt];
  \end{scope}
 \end{tikzpicture}
 \caption{Construction of $H$: before gluing (left) and after gluing (right).}
 \label{fig:handle-gluing}
\end{figure}

Consider a right-angled regular hyperbolic octagon with sides $x,u,v,w,y,\bar{u},\bar{v},\bar{w}$ in counterclockwise order, and oriented in that direction. Pair those sides with orientation-reversing isometries according to the identification
\begin{equation*}
    u\leftrightarrow\bar u,\qquad v\leftrightarrow\bar v,\qquad w\leftrightarrow\bar w,
\end{equation*}
and denote by $H$ the resulting quotient, which is homeomorphic to a genus 1 surface with two boundary components; see \Cref{fig:handle-gluing}.

Write $e\coloneqq\alpha_0$ in $Q_0$. Attach $H$ to $Q_0$ via orientation-preserving isometries identifying $x\leftrightarrow e^{-1}$ and $y\leftrightarrow e$, to obtain the octagon complex
\begin{equation*}
    Q_*\coloneqq Q_0\cup_e H
\end{equation*}
with one extra octagon face of boundary word $e^{-1}uvweu^{-1}v^{-1}w^{-1}$. Let $f_0\colon Q_0\monic Q_*$ be the inclusion.

Fix $k\geq 1$ and let $Y'_0$ be obtained from $Y_0$ by cutting along $\alpha_0$, i.e.\ undoing the identification $\alpha_0\leftrightarrow\bar\alpha_0$. Attach $k$ copies $H_1,\dots, H_k$ of $H$, with sides labelled with corresponding subscripts, to $Y'_0$ under the following identifications
\begin{equation*}
    \alpha_0\leftrightarrow x_1,\qquad y_j\leftrightarrow x_{j+1}\quad (1\leq j<k),\qquad y_k\leftrightarrow\bar\alpha_0,
\end{equation*}
via orientation-reversing isometries. The resulting octagon complex $Y_k$ is a closed oriented surface of genus $431+k$. Let $\Gamma_{431+k}\coloneqq\pi_1(Y_k)$.

Let $Q_k\coloneqq Z\cup_R Y_k$ with gluing embeddings $R\monic Z$ and $R\monic Y_k$ (defined as for $R\monic Y_0$, sending $R$ to the $\delta$-edges). The Seifert--van Kampen theorem yields $K_k\coloneqq\pi_1(Q_k)$ as the pushout
\[
    \begin{tikzcd}
        R \arrow[r, hook] \arrow[d, hook] & Z \arrow[d, hook] \\
        Y_k \arrow[r, hook] & Q_k\coloneqq Z\cup_R Y_k,
    \end{tikzcd}
    \qquad\qquad
    \begin{tikzcd}
        F \arrow[r] \arrow[d] & G \arrow[d] \\
        \Gamma_{431+k} \arrow[r] & K_k=G\ast_F \Gamma_{431+k}.
    \end{tikzcd}
\]

Consider the following map of octagon complexes
\begin{equation*}
    f_k\colon Q_k\to Q_*
\end{equation*}
induced by the inclusion $Z\to Q_0$ and by the map $Y_k\to Q_*$ sending the complement of the new octagons to their copies in $Q_0$, and sending each $H_j$ isometrically to $H$ (which sends the edges $x_{j}^{-1},y_{j},\alpha_0,\bar\alpha_0^{-1}$ to $e$, isometrically and preserving the orientation). Both maps coincide on $R$, hence $f_k$ is well-defined.

\begin{proposition}\label{prop:Qk_higher_genus}
    Let $k\geq 1$. The following assertions hold:
    \begin{enumerate}[label=(\roman*)]
        \item The complex $Q_*$ is a one-vertex octagon complex with a simple and triangle-free vertex link, with immersed face boundaries and all edges contained in boundaries. Thus, it is a locally \cat{-1}\ octagon complex and $K_*\coloneqq\pi_1(Q_*)$ is torsion-free and word-hyperbolic.
        \item The embeddings (or markings) $R\monic Z$ and $R\monic Y_k$ induce injective homomorphisms between fundamental groups. Thus $K_k=G\ast_F \Gamma_{431+k}$ is an amalgamated free product and $G$, $\Gamma_{431+k}$ embed in $K_k$. Also, the image of the injection $F\to\Gamma_{431+k}$ is quasi-convex.
        \item The map $f_k$ is a local isometry and the induced homomorphism $(f_k)_*\colon K_k\to K_*$ is injective with quasi-convex image. The same holds for $f_0$.
    \end{enumerate}
\end{proposition}
\begin{proof}
    In this proof, let us write the half-edges as $a^{\pm}\coloneqq [a]_{\pm}$. The link of $Q_*$ is the union of the link of $Q_0$ with the two new paths
    \begin{align}
        P_x\coloneqq e^+\text{--}u^+\text{--}v^-\text{--}w^+\text{--}e^-,\label{eq:Qstar_from_Q0_first}\\
        P_y\coloneqq e^+\text{--}w^-\text{--}v^+\text{--}u^-\text{--}e^-.\label{eq:Qstar_from_Q0_second}
    \end{align}
    Their internal vertices are new and distinct, so this does not create any loop, parallel edge or triangle. See \Cref{fig:vertex_link_Qstar}. The new face boundary is immersed, and every edge belongs to a face. Thus (i) follows from \Cref{prop:construction_of_Q_0,prop:hyperbolicity_of_octagon_complex}.

    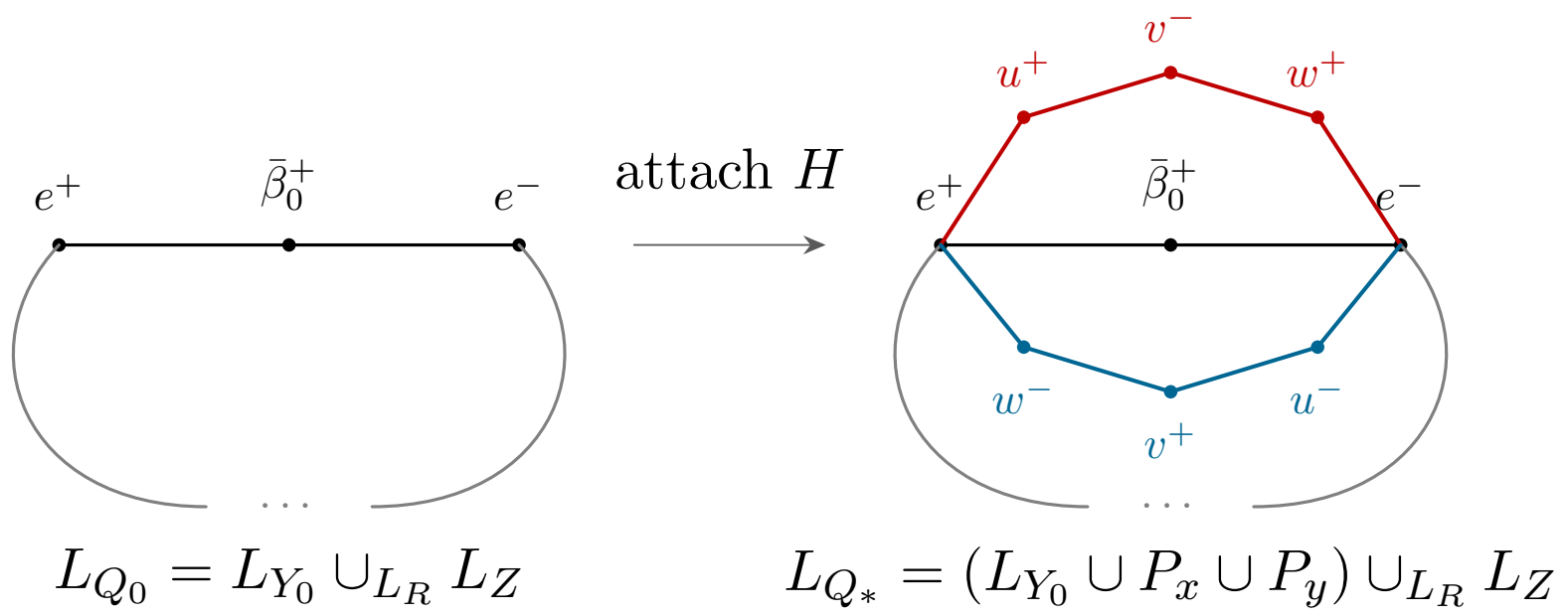
\begin{figure}[ht]
    \expandafter\def\csname @captype\endcsname{figure}
    \centering
    \scalebox{0.85}{%
    \begin{tikzpicture}[line cap=round,line join=round,
        every node/.style={font=\small},
        old link/.style={line width=.7pt},
        x link/.style={red!75!black,line width=.9pt},
        y link/.style={MidnightBlue,line width=.9pt},
        omitted link/.style={black!50,line width=.7pt}]
      \foreach \shift in {0,6.9} {
        \begin{scope}[xshift=\shift cm]
          \draw[old link] (-1.8,0)--(0,0)--(1.8,0);
          \foreach \a/\lab in {-1.8/{e^+},0/{\bar\beta_0^+},1.8/{e^-}} {
            \fill (\a,0) circle[radius=1.5pt];
            \node[above=4pt] at (\a,0) {$\lab$};
          }
          \draw[omitted link] (-1.8,0)
            .. controls (-2.55,-.8) and (-2.1,-2.05) .. (-.65,-2.05);
          \draw[omitted link] (1.8,0)
            .. controls (2.55,-.8) and (2.1,-2.05) .. (.65,-2.05);
          \node[text=black!50] at (0,-2.05) {$\cdots$};
        \end{scope}
      }
      \begin{scope}[xshift=6.9cm]
        \draw[x link] (-1.8,0)--(-1.15,1)--(0,1.35)--(1.15,1)--(1.8,0);
        \foreach \a/\b/\lab in {-1.15/1/{u^+},0/1.35/{v^-},1.15/1/{w^+}} {
          \fill[red!75!black] (\a,\b) circle[radius=1.5pt];
          \node[above=3pt,text=red!75!black] at (\a,\b) {$\lab$};
        }
        \draw[y link] (-1.8,0)--(-1.15,-.8)--(0,-1.15)--(1.15,-.8)--(1.8,0);
        \foreach \a/\b/\lab in {-1.15/-.8/{w^-},0/-1.15/{v^+},1.15/-.8/{u^-}} {
          \fill[MidnightBlue] (\a,\b) circle[radius=1.5pt];
          \node[below=3pt,text=MidnightBlue] at (\a,\b) {$\lab$};
        }
      \end{scope}
      \draw[-{Stealth[length=5pt,width=4.2pt]},black!65] (2.7,0)--(4.2,0);
      \node[font=\normalsize,scale=1.176,align=center] at (3.45,.6) {attach $H$};
      \node[font=\normalsize,scale=1.176] at (0,-2.6) {$L_{Q_0}=L_{Y_0}\cup_{L_R}L_Z$};
      \node[font=\normalsize,scale=1.176] at (6.9,-2.6) {$L_{Q_*}=(L_{Y_0}\cup P_x\cup P_y)\cup_{L_R}L_Z$};
    \end{tikzpicture}
    }
    \caption{Left: The vertex link $L_{Q_0}$ of $Q_0$. Right: The vertex link $L_{Q_*}$ of $Q_*$. The grey curve with dots represents the remaining part of the vertex links, which coincide in $Q_0$ and $Q_*$: it consists of the rest of $L_{Y_0}$, to which $L_Z$ is attached along $L_R$.}
    \label{fig:vertex_link_Qstar}
    \end{figure}

    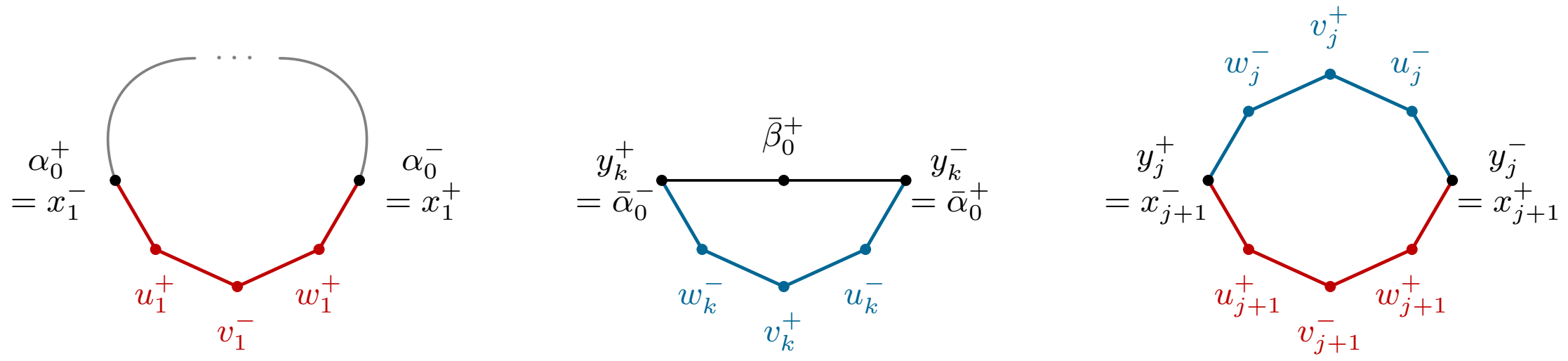
\begin{figure}[ht!]
    \centering
    \begin{tikzpicture}[line cap=round,line join=round,
        every node/.style={font=\footnotesize},
        old link/.style={line width=.7pt},
        x link/.style={red!75!black,line width=.9pt},
        y link/.style={MidnightBlue,line width=.9pt},
        omitted link/.style={black!50,line width=.7pt}]
      \begin{scope}
        \draw[omitted link] (0,0)
          .. controls (-.24,.7) and (.17,1.15) .. (.75,1.15);
        \draw[omitted link] (2.3,0)
          .. controls (2.54,.7) and (2.13,1.15) .. (1.55,1.15);
        \node[text=black!50] at (1.15,1.15) {$\cdots$};
        \draw[x link] (0,0)--(.38,-.65)--(1.15,-1)--(1.92,-.65)--(2.3,0);
        \foreach \a/\b in {0/0,2.3/0} {\fill (\a,\b) circle[radius=1.5pt];}
        \foreach \a/\b in {.38/-.65,1.15/-1,1.92/-.65}
          {\fill[red!75!black] (\a,\b) circle[radius=1.5pt];}
        \node[left=3pt,align=center] at (0,0) {$\alpha_0^+$\\$=x_1^-$};
        \node[right=3pt,align=center] at (2.3,0) {$\alpha_0^-$\\$=x_1^+$};
        \foreach \a/\b/\lab in {.38/-.65/{u_1^+},1.15/-1/{v_1^-},1.92/-.65/{w_1^+}}
          {\node[below=3pt,text=red!75!black] at (\a,\b) {$\lab$};}
      \end{scope}
      \begin{scope}[xshift=5.15cm]
        \draw[old link] (0,0)--(1.15,0)--(2.3,0);
        \draw[y link] (0,0)--(.38,-.65)--(1.15,-1)--(1.92,-.65)--(2.3,0);
        \foreach \a/\b in {0/0,1.15/0,2.3/0}
          {\fill (\a,\b) circle[radius=1.5pt];}
        \foreach \a/\b in {.38/-.65,1.15/-1,1.92/-.65}
          {\fill[MidnightBlue] (\a,\b) circle[radius=1.5pt];}
        \node[above=3pt] at (1.15,0) {$\bar\beta_0^+$};
        \node[left=3pt,align=center] at (0.2,0) {$y_k^+$\\$=\bar\alpha_0^-$};
        \node[right=3pt,align=center] at (2.1,0) {$y_k^-$\\$=\bar\alpha_0^+$};
        \foreach \a/\b/\lab in {.38/-.65/{w_k^-},1.15/-1/{v_k^+},1.92/-.65/{u_k^-}}
          {\node[below=3pt,text=MidnightBlue] at (\a,\b) {$\lab$};}
      \end{scope}
      \begin{scope}[xshift=10.3cm]
        \draw[y link] (0,0)--(.38,.65)--(1.15,1)--(1.92,.65)--(2.3,0);
        \draw[x link] (0,0)--(.38,-.65)--(1.15,-1)--(1.92,-.65)--(2.3,0);
        \foreach \a/\b in {0/0,2.3/0} {\fill (\a,\b) circle[radius=1.5pt];}
        \foreach \a/\b in {.38/.65,1.15/1,1.92/.65}
          {\fill[MidnightBlue] (\a,\b) circle[radius=1.5pt];}
        \foreach \a/\b in {.38/-.65,1.15/-1,1.92/-.65}
          {\fill[red!75!black] (\a,\b) circle[radius=1.5pt];}
        \node[left=3pt,align=center] at (0.25,0) {$y_j^+$\\$=x_{j+1}^-$};
        \node[right=3pt,align=center] at (2.1,0) {$y_j^-$\\$=x_{j+1}^+$};
        \foreach \a/\b/\lab in {.38/.65/{w_j^-},1.15/1/{v_j^+},1.92/.65/{u_j^-}}
          {\node[above=3pt,text=MidnightBlue] at (\a,\b) {$\lab$};}
        \foreach \a/\b/\lab in {.38/-.65/{u_{j+1}^+},1.15/-1/{v_{j+1}^-},1.92/-.65/{w_{j+1}^+}}
          {\node[below=3pt,text=red!75!black] at (\a,\b) {$\lab$};}
      \end{scope}
    \end{tikzpicture}
    \caption{The three types of vertex links of $Y_k$ and $Q_k$. From left to right: the main vertex link, the cycle of length $6$, and the $k-1$ cycles of length $8$.}
    \label{fig:vertex_links_of_Qk}
    \end{figure}

    In $Y_k$, all marked half-edges (i.e., in the image of $R\to Y_k$) remain incident to the main vertex. Its main vertex link is obtained from that of $Y_0$ (see \eqref{eq:link_vertex_Y0}) by replacing the two-edge path
    \begin{equation*}
        \alpha_0^+\text{--}\bar\beta_0^+\text{--}\alpha_0^-
    \end{equation*}
    with the four-edge path
    \begin{equation}\label{eq:main_link_vertex_Yk}
        \alpha_0^+=x_1^-\text{--}u_1^+\text{--}v_1^-\text{--}w_1^+\text{--}x_1^+=\alpha_0^-.
    \end{equation}
    The other vertex links in $Y_k$ are one cycle of length $6$ and $k-1$ cycles of length $8$. Hence all links are simple and triangle-free. Moreover, distances between marked half-edges do not decrease in $Y_k$ compared to their distances in $Y_0$, and no new two-edge path joins marked half-edges. Thus the two conditions of \Cref{lem:attachment_criterion}, already verified for $Y_0$, still hold for $Y_k$. Face boundaries remain immersed and every edge lies in a face. Since $R\monic Z$ is unchanged, \Cref{lem:attachment_criterion} applies and proves (ii) as in \Cref{prop:construction_of_Q_0}.

    Finally, $Q_k$ has $k+1$ vertices: one ``main'' vertex, one vertex whose link is a cycle of length $6$ and $k-1$ vertices whose links are cycles of length $8$; see \Cref{fig:vertex_links_of_Qk}. The new vertices originate from $Y_k$. The main vertex link is obtained from the vertex link of $Q_0$ by deleting $\bar\beta_0^+$ and adjoining the path $P_x$ from \eqref{eq:Qstar_from_Q0_first}. The map $f_k$ sends $u_j,v_j,w_j$ to $u,v,w$, respectively. Looking at \Cref{fig:vertex_link_Qstar,fig:vertex_links_of_Qk}, we see that $f_k$ induces injective maps onto full subgraphs between the links. This also holds for $k=0$, as the inclusion of $L_{Q_0}$ into $L_{Q_*}$ can be visualized in \Cref{fig:vertex_link_Qstar}. Together with the facewise isometries and local injectivity along edges, this proves that $f_k$ is a local isometry. By \Cref{prop:local_isometries_are_pi_1_injective}, it induces an injection and lifts to an isometric embedding of universal covers. Its image is convex, and cocompactness of the $K_k$-action on this image implies that $(f_k)_*(K_k)$ is quasi-convex in $K_*$. This proves (iii).
\end{proof}

\begin{proof}[Proof of \Cref{thm:all_genera}]
    By \Cref{prop:Qk_higher_genus}(i), the universal covering $\widetilde{Q_*}$ endowed with the metric where each of its faces is a regular Euclidean octagon satisfies the hypotheses of \Cref{thm:Desgroseilliers--Haglund}, as in the proof of \Cref{thm:main}. This yields $q_*\in\N$ and a discrete faithful convex cocompact representation $\sigma_*\colon K_*\to\Isom(\HH^{q_*})$. It is purely loxodromic since $K_*$ is torsion-free, and the translation lengths of its non-trivial elements are bounded below by a positive constant.

    Let $g=431+k$ with $k\geq 0$. By \Cref{prop:Qk_higher_genus}(ii,iii) (and \Cref{prop:construction_of_Q_0} for $k=0$), $(f_k)_*\restr{\Gamma_g}\colon\Gamma_g\to K_*$ is injective; denote its image by $\Gamma'_g$. Thus $\rho_g\coloneqq \sigma_*\restr{\Gamma'_g}$ is discrete, faithful and purely loxodromic, with the same lower bound on translation lengths.

    Similarly, $\sigma_*\restr{(f_k)_*G}$ is discrete, faithful and purely loxodromic, hence $\sigma_*\restr{(f_k)_*F}$ has no \ct\ map by \Cref{prop:explicit_obstruction}. Since $F\monic\Gamma_g$ is quasi-convex, the same holds for $\rho_g$. In particular, $\rho_g$ is not convex cocompact.
\end{proof}

\section{Ray Cannon--Thurston maps}\label{sec:ray_CT_map}

We now describe more precisely the boundary behaviour of the examples constructed above. Every geodesic ray has an image converging to a point at infinity, even though the orbit map does not extend continuously to the boundary. This is captured by the notion of a ray Cannon--Thurston map, introduced by Halder, Mahan Mj and Sardar~\cite{HMS2025rayCTmap}. We recall their definition and prove that it applies to our examples.

\begin{definition}[\!\!{\cite[Definition 2.6]{HMS2025rayCTmap}}]
    Let $i\colon Y\to X$ be a map between Gromov hyperbolic geodesic metric spaces. A (quasi-)geodesic ray $r\colon [0,\infty)\to Y$ is called a \emph{landing ray} for $i$ if there exists $\zeta\in\partial X$ such that
    \begin{equation*}
        i(r(t))\longrightarrow\zeta\qquad\text{as }t\to\infty.
    \end{equation*}
    In that case, we say that $i\circ r$ lands at $\zeta$ in $\partial X$.

    If every (quasi-)geodesic ray in $Y$ is a landing ray, then this defines a map
    \begin{equation*}
        \partial_r i\colon\partial Y\to\partial X\colon \eta\mapsto\lim_{t\to\infty}i(r_{\eta}(t))
    \end{equation*}
    where $r_{\eta}\colon [0,\infty)\to Y$ is a quasi-geodesic ray ending at $\eta$. The map $\partial_r i$ is then called a \emph{ray \ct\ map} for the map $i\colon Y\to X$.

    An inclusion of hyperbolic groups $H\leq G$ admits a ray \ct\ map if the inclusion $H\to G$, for any choice of word metrics with respect to finite generating sets, does. A representation $\rho\colon \Gamma\to\Isom(X)$ of a hyperbolic group $\Gamma$ on a Gromov hyperbolic space $X$ admits a ray \ct\ map if one of its orbit maps does.
\end{definition}

\begin{proposition}
    The inclusions $F\leq G$ and $\Gamma_{431+k}\to K_k$, for $k\geq 0$, all admit ray \ct\ maps.
\end{proposition}
\begin{proof}
    First, let us prove that $F\leq G$ admits a ray \ct\ map. This is a direct consequence of the criterion of Halder, Mahan Mj and Sardar~\cite[Theorem 1.3]{HMS2025rayCTmap}, applied to the following choices of subgroups:
    \begin{itemize}
        \item The subgroups $A\coloneqq\langle c_i,d_j\rangle$, $B\coloneqq \langle b,c_i\rangle$, $C=\langle c_i\rangle$ of $G$;
        \item The subgroups $U\coloneqq\langle d_j\rangle$, $V\coloneqq \langle b\rangle$, $W\coloneqq\{1\}$ of $F$.
    \end{itemize}
    Indeed, we have $G=A\ast_{C}B$, $F=U\ast_{W}V$, and the only conditions that require checking are that $A$ and $B$ are hyperbolic, that $C\leq A$ is malnormal and quasi-convex and that both inclusions $U\leq G$ and $V\leq G$ admit a ray \ct\ map. The rest of the conditions are satisfied: $U/W\to A/C$ is injective because $U\cap C=\{1\}$ (\Cref{lem:double_HNN}), $V/W\to B/C$ is injective because $\langle b\rangle\cap C=\{1\}$ and $W=\{1\}=U\cap V$.

    Hyperbolicity of $A$ and $B$ follows from realizing them as octagon subcomplexes of $Z$ (see \Cref{prop:combinatorial_description_Z}) and verifying that their vertex links satisfy \Cref{prop:hyperbolicity_of_octagon_complex}.

    The inclusion $C\leq A$ is given explicitly by the iterated HNN extension
    \begin{equation*}
        \langle c_i\rangle\to\langle c_i,d_j\mid c_i^{-1}d_jc_i=\psi_i(d_j)\rangle,
    \end{equation*}
    where $C=\langle c_i\rangle$ is freely generated by the letters $c_i$. Thus, by \cite[Corollary 6.6]{HMS2025rayCTmap}, $C\leq A$ is indeed malnormal and quasi-convex.

    The inclusion $\langle b\rangle=V\leq G$ is quasi-convex, hence it admits a ray \ct\ map. The inclusion $U\leq G$ is the composition $U\leq A\leq G$, where both $U\leq A$ and $A\leq G$ actually admit \ct\ maps: the former by \cite[Lemma 6.8(1)]{HMS2025rayCTmap} and the latter by \cite[Proof of Theorem 6.7]{HMS2025rayCTmap}; both also follow from \cite[Corollary 3.11]{MR1622603}.

    Therefore, by \cite[Theorem 1.3]{HMS2025rayCTmap}, the inclusion $F\leq G$ admits a ray \ct\ map.

    Let $k\geq 0$ and write $g\coloneqq 431+k$. Let us show that $\Gamma_{g}\to K_k$ admits a ray \ct\ map, using \Cref{prop:ray_CT_map_stable_under_malnormal_quasi-convex_extensions}. It suffices to check that $F\to\Gamma_{g}$ is (almost) malnormal and quasi-convex. Quasi-convexity holds by \Cref{prop:construction_of_Q_0,prop:Qk_higher_genus}(ii). For the same reason, the embedding $R\monic Y_k$ lifts to an isometric embedding $\widetilde{R}\monic\widetilde{Y_k}$, realizing $\widetilde{R}$ as a convex embedded tree $T$ in $\widetilde{Y_k}$. Let $\gamma\in\Gamma_g-F$, and suppose that $1\neq h\in F\cap \gamma F\gamma^{-1}$. Then the action of $h$ on $\widetilde{Y_k}$ preserves both trees $T$ and $\gamma T$, and acts loxodromically on each one, say along translation axes $\ell$ and $\ell'$. But we must have $\ell=\ell'$, as $h$ acts loxodromically on the \cat{-1}\ space $\widetilde{Y_k}$, contradicting the fact that $T$ and $\gamma T$ are disjoint: indeed, they are distinct connected components of the preimage of $R$ in $\widetilde{Y_k}$. This proves that $F\leq\Gamma_g$ is malnormal. Thus, by \Cref{prop:ray_CT_map_stable_under_malnormal_quasi-convex_extensions}, and since $F\leq G$ admits a ray \ct\ map, we conclude that $\Gamma_g\leq K_k$ admits a ray \ct\ map.
\end{proof}

\begin{proof}[Proof of \Cref{prop:ray_CT_maps}]
    Finally, $K_*$ acts convex cocompactly on $\HH^{q_*}$ by the proof of \Cref{thm:all_genera} and by \Cref{thm:Desgroseilliers--Haglund}, hence there is a \ct\ map $\partial K_*\monic \partial\HH^{q_*}$. Each $K_k\to K_*$ is quasi-convex by \Cref{prop:Qk_higher_genus}(iii), hence the orbit maps of the representations $F,\Gamma_{431+k}\to\Isom(\HH^{q_*})$ all admit ray \ct\ maps, by the above Proposition and postcomposing with the embedding $\partial K_k\monic\partial K_*$. This also holds for $\sigma_0\restr{\Gamma_{431}}$.
\end{proof}

\appendix

\section{An extension property for ray \ct\ maps}\label{sec:extension_ray_CT_map}
In~\cite[Theorem A.1]{HMS2025rayCTmap}, it is shown that having a \ct\ map is stable under almost malnormal quasi-convex hyperbolic extension. The following proposition establishes the analogous statement for ray \ct\ maps.

We recall that a subgroup $F\leq\Gamma$ is \emph{almost malnormal} if, for every $\gamma\in\Gamma-F$, $F\cap \gamma F\gamma^{-1}$ is finite.

\begin{proposition}\label{prop:ray_CT_map_stable_under_malnormal_quasi-convex_extensions}
    Suppose $K=\Gamma\ast_{F} G$ is a free product with amalgamation such that $\Gamma,G$ are hyperbolic and $F$ is infinite, almost malnormal and quasi-convex in $\Gamma$. Then $\Gamma\leq K$ admits a ray \ct\ map if and only if $F\leq G$ does.
\end{proposition}
\begin{proof}
    Since $F$ is almost malnormal and quasi-convex in $\Gamma$, \cite[Lemma 2.18]{HMS2025rayCTmap} implies that $K$ is hyperbolic and that $G$ is almost malnormal and quasi-convex in $K$; see also~\cite{MR1152226}. Moreover, quasi-convexity of $F\leq\Gamma$ implies that $F$ is hyperbolic and that the inclusion $F\to\Gamma$ is a quasi-isometric embedding. For this proof, let us denote the different group inclusions as in the diagram
    \[
        \begin{tikzcd}
            F \arrow[r, hook, "\iota"] \arrow[d, hook, "\kappa"] & G \arrow[d, hook, "\kappa'"] \\
            \Gamma \arrow[r, hook, "\iota'"] & K=G\ast_F \Gamma.
        \end{tikzcd}
    \]

    Suppose that $\Gamma\leq K$ admits a ray \ct\ map. Let $\eta\in\partial F$ and consider the ray $r$ from $1$ to $\eta$ in $F$. Since $F\leq\Gamma$ is quasi-convex, the ray $\kappa\circ r$ is a quasi-geodesic ray in $\Gamma$. Thus, it lands at a point $\xi\in \partial K$ by hypothesis, i.e.
    \begin{equation*}
        \kappa'\circ\iota\circ r(t)=\iota'\circ\kappa\circ r(t)\longrightarrow\xi.
    \end{equation*}
    Since $\kappa'\colon G\monic K$ is quasi-convex, $\xi$ must lie in $\partial G$ and we have $\iota\circ r(t)\longrightarrow \xi$. This shows that $F\leq G$ admits a ray \ct\ map.

    Now, suppose that $F\leq G$ admits a ray \ct\ map. We shall prove that $\Gamma\leq K$ does too. The rest of the argument adapts \cite[Proof of Theorem 3.1, Step 3, Case 2]{HMS2025rayCTmap}.
    
    Let us introduce some notation, following \cite[Proof of Theorem A.1]{HMS2025rayCTmap}; see \cite{MR1650094,MR2780738,MR3646028} for more details. Let $\widehat{K}$ denote the electrification of the Cayley graph of $K$ obtained by coning off the cosets of $G$ in $K$. More precisely, for each coset $kG$ in $K/G$, a vertex $v_{kG}$ is added and joined to every vertex of $kG\subset K$ by an edge of length $1$. Similarly, let $\widehat{\Gamma}$ denote the electrification of the Cayley graph of $\Gamma$ by coning off the cosets of $F$ in $\Gamma$. As in \cite[Proof of Theorem A.1]{HMS2025rayCTmap}, the induced map $\widehat{\iota}'\colon\widehat{\Gamma}\to\widehat{K}$ is a quasi-isometric embedding. Quasi-convexity allows us to replace each passage through a cone vertex by a uniform quasi-geodesic contained in the corresponding coset; we call these the replacement intervals. Since $F\leq\Gamma$ and $G\leq K$ are almost malnormal and quasi-convex, \cite[Theorem 7.11]{MR2922380} implies that $\Gamma$ and $K$ are strongly hyperbolic relative to $F$ and $G$, respectively. Consequently, de-electrifications of electric quasi-geodesics without back-tracking are uniform quasi-geodesics \cite[Lemma 2.15]{MR3646028}.

    Let $\eta\in\partial\Gamma$, and let $r$ be a ray from $1$ to $\eta$ in $\Gamma$. Suppose that $\eta\in \gamma\partial F$ for some $\gamma\in\Gamma$. Since $\gamma F$ is quasi-convex in $\Gamma$, the ray $r$ stays at bounded Hausdorff distance from a ray $\gamma\cdot \kappa\circ r'$ where $r'$ is a ray in $F$. Since $F\leq G$ admits a ray \ct\ map, the ray $\iota\circ r'$ lands at a point $\xi\in\partial G$. By quasi-convexity of $G\leq K$, the ray $\kappa'\circ\iota\circ r'$ lands at $\xi$ in $\partial K$ too. By $\Gamma$-equivariance of $\iota'$, this gives
    \begin{equation*}
        \iota'\bigl(\gamma\cdot \kappa (r'(t))\bigr)\longrightarrow\gamma\cdot\xi\in\gamma\cdot\partial G\subset\partial K.
    \end{equation*}
    Because $r$ is at bounded distance from $\gamma\cdot(\kappa\circ r')$ and $\iota'$ is Lipschitz, we get $\iota'\circ r(t)\longrightarrow \gamma\cdot\xi$ too.

    Now, suppose that $\eta$ does not lie in the subset $\bigcup_{\gamma\in\Gamma}\gamma\partial F$ of $\partial\Gamma$. By \cite[Proposition 3.14]{MR4844109} and local finiteness of the family of left cosets of $F$ \cite[Example 3.13]{MR4844109}, the ray $\widehat{r}$, obtained by postcomposing $r$ with $\Gamma\to\widehat{\Gamma}$, has unbounded image in $\widehat{\Gamma}$. By \cite[Lemma 3.9]{MR4844109}, there is an electric quasi-geodesic ray $\widehat{\alpha}$ in $\widehat{\Gamma}$ at bounded Hausdorff distance from $\widehat{r}$. We can assume that $\widehat{\alpha}$ is without back-tracking (i.e.\ it does not enter any cone vertex twice).

    Let $\alpha$ be the de-electrification of $\widehat{\alpha}$ in $\Gamma$ \cite[Definition 3.3]{MR4844109}. It is a quasi-geodesic in $\Gamma$ \cite[Lemma 2.15]{MR3646028} and, by \cite[Theorem 3.11]{MR4844109}, $\alpha(t)\longrightarrow \eta$, hence $\alpha$ and $r$ lie at bounded Hausdorff distance from one another.

    Since $\Gamma\cap \gamma G=\gamma F$ for each $\gamma\in\Gamma$, the quasi-geodesic ray $\widehat{\beta}\coloneqq\widehat{\iota}'\circ\widehat{\alpha}$ in $\widehat{K}$ remains without back-tracking, hence its de-electrification $\beta$ in $K$ is a quasi-geodesic ray, landing at some point $\xi\in\partial K$.

    The rays $\iota'\circ\alpha$ and $\beta$ coincide outside their replacement intervals. If $\widehat{\beta}$ traverses only finitely many cone vertices, the two paths have a common tail so $\iota'\circ\alpha$ lands at $\xi$ too. Otherwise, let $C_i=p_iG$ be the cosets whose cone vertices $v_{C_i}$ are traversed by $\widehat{\beta}$, in order, where $p_i\in\Gamma$ is the entry point of the corresponding replacement interval. Quasi-geodesicity of $\widehat{\beta}$ gives (using that $K\to\widehat{K}$ is Lipschitz for the first inequality)
    \begin{equation*}
        d_K(1,C_i)\geq d_{\widehat{K}}(1,C_i)\geq d_{\widehat{K}}(1,v_{C_i})-1\longrightarrow +\infty.
    \end{equation*}
    Also, $p_i\longrightarrow\xi$ along the quasi-geodesic ray $\beta$. Since $G$ is quasi-convex in $K$, its cosets $C_i$ are uniformly quasi-convex. Thus \cite[Lemma 2.13]{HMS2025rayCTmap} implies that $C_i\longrightarrow\xi$. Each corresponding replacement interval of $\iota'\circ\alpha$ lies in $p_iF\subset p_iG=C_i$, so $\iota'\circ\alpha$ also lands at $\xi$. Finally, $\alpha$ and $r$ lie at bounded Hausdorff distance in $\Gamma$ and $\iota'$ is Lipschitz, hence $\iota'\circ r$ lands at $\xi$ as well. This shows that $\Gamma\leq K$ has a ray \ct\ map.
\end{proof}

\bibliographystyle{alpha}
\bibliography{bib}
\end{document}